\documentclass[1p]{ELSE-article}
\usepackage[cp1251]{inputenc}
\usepackage[russian,english]{babel}
\usepackage{amsmath,amssymb,graphicx,color}
\usepackage[colorlinks,pdfauthor={Yu. Brezhnev},
            pdfwindowui=false,%pdfmenubar=false,%pdftoolbar=false,
            bookmarksopen=false,bookmarks=false]{hyperref}
\def~{\unskip\nobreak\hspace{0.27778em}\ignorespaces}
\def\,{\ifmmode\mskip+1.5mu\else\kern+.08333em\fi\relax}
\def\!{\ifmmode\mskip-1.5mu\else\kern-.08333em\fi\relax}
\newdimen{\FontSize} \delimiterfactor=990
\def\cdot{\,{\mathchar"2201}\,}

\makeatletter
\def\dot#1{\mathbf{\mathaccent"705F\mathnormal{#1}}}

\def\dddot#1{{\mathop{{}#1}\limits^{\vbox to-1.5\ex@{\kern-\tw@\ex@
\hbox{\larger[2]\rm.\kern-0.12em.\kern-0.12em.}\vss}}}} \makeatother

\def\ds{\displaystyle}

\def\s{\scriptscriptstyle}
\def\ie{i.\,\,e.}

\def\eg{\text{e.\,\,g.}}
\def\re{{\,\mathrm{e}}}
\def\ri{{\mathrm{i}}}

\def\Hp{\mathbb{H}^+}

\def\ode{{\smaller[2]\text{ODE}}} % работает в \italic

\def\wpp {\wp\Smaller{'}}

\def\omp {\omega\Smaller{'}}

\def\ep {e\Smaller{'}}
\def\epp{e\Smaller{'\mskip-1.2mu'}}

\def\g#1{{g_{#1}^{}}}

\def\ellK{{\mathsf{K}}}

\def\II{{\textbf{\textsc{I\kern -0.15em I}}}}
\def\III{{\textbf{\textsc{I\kern -0.17em I\kern -0.17em I}}}}
\def\IIIab{\III\!\raise0.2ex\hbox{$\begin{smallmatrix}
\alpha\\[-0.1ex]\varrho\end{smallmatrix}$}\kern-0.1em}

\def\DEF{\mathrel{\vcenter{\hbox{$:$}}{=}}}
\def\FED{\mathrel{{=}\vcenter{\hbox{$:$}}}}
\def\=={\mathrel{\phantom{=}}}

\def\sm{\mathrel{\vcenter{\hbox{\scalebox{0.8}[1]{$\scriptstyle-$}}}}}
\def\hence{\scalebox{1.5}[1]{\hbox{$\Rightarrow$}}}
\def\?{\textrm{\protect\footnotesize$\red\mathchar"446$}}

\def\vphi{\smash[b]{\raise-0.22ex\hbox{\Smaller[2]{\s\boldsymbol-}}
  \mkern-8.00mu\raise0.435ex\hbox{\scalebox{1}[0.9]{$\varphi$}}}\relax}
\def\vpsi{\smash[b]{\raise-0.22ex\hbox{\Smaller[2]{\s\boldsymbol-}}
  \mkern-8.95mu\raise0.408ex\hbox{\scalebox{1}[0.9]{$\psi   $}}}\relax}
\def\russian{\selectlanguage{russian}}
\def\english{\selectlanguage{english}}
\def\END{\end{document}}

\def\Mfrac#1#2{\def\FRAC{\hbox{\smaller[1]$\ds\frac{#1}{#2}$}}
 \ifdim\FontSize=10pt\raise0.057ex\FRAC\else
{\ifdim\FontSize=11pt\raise0.050ex\FRAC\else\raise0.051ex\FRAC\fi}\fi}

\def\mfrac#1#2{\def\FRAC{\hbox{\smaller[2]$\ds\frac{#1}{#2}$}}
 \ifdim\FontSize=10pt\raise0.1166ex\FRAC\else
{\ifdim\FontSize=11pt\raise0.104ex\FRAC\else\raise0.0968ex\FRAC\fi}\fi}

\newcommand{\AB}[2]{\nobreak\mskip-1mu
\raise0.12ex\hbox{\smaller[2]\renewcommand{\arraystretch}{0.42}%
\setlength{\arraycolsep}{0pt}\raise0.16ex\hbox{$\big[$}%
$\begin{array}{c}\scriptstyle#1\\ \scriptstyle#2\end{array}$%
\raise0.16ex\hbox{$\big]$}}\mskip-1mu\relax}

\newcommand{\mbig}[2][4]{\vcenter{\hbox{%
\ifcase#1$#2$\or\small$\big#2$\or$\big#2$\or\large$\big#2$%
\or\footnotesize$\Big#2$\or\small$\Big#2$\or$\Big#2$\or\large$\Big#2$%
\or$\bigg#2$\or\large$\bigg#2$\or$\Bigg#2$\or\large$\Bigg#2$%
\or\Large$\Bigg#2$\else\huge$\bigg#2$\fi}}\relax}

\newdimen{\upp} \newdimen{\Dim} \newdimen{\DDim}
\newcommand{\pow}[4][\displaystyle]{
\settowidth{\upp}{\hbox{$#1 {}^{#4}$}}%
\settowidth{\Dim}{\hbox{$#1{#2}$}}%
\settowidth{\DDim}{\hbox{$#1{#2}_{#3}$}}%
\addtolength{\DDim}{-\Dim} {#1{#2}_{#3}\kern-\DDim {}^{#4}}
\addtolength{\DDim}{-\upp}\ifdim\DDim>0pt \kern+\DDim \fi\relax}

\makeatletter % \@tempcnta и \@currsize - внутренние LaTeX's переменные
\providecommand{\larger}[1][1]{% работает в точности как \larger в amsart
\ifmmode\relax\else\@tempcnta\ifx\@currsize\normalsize 4\else
 \ifx\@currsize\small 3\else\ifx\@currsize\footnotesize 2\else
 \ifx\@currsize\large 5\else\ifx\@currsize\Large 6\else
 \ifx\@currsize\LARGE 7\else\ifx\@currsize\scriptsize 1\else
 \ifx\@currsize\tiny  0\else\ifx\@currsize\huge 8\else
 \ifx\@currsize\Huge  9\else 4 \fi\fi\fi\fi\fi\fi\fi\fi\fi\fi
\advance\@tempcnta#1\relax\ifnum\@tempcnta<\z@ \@tempcnta\z@ \fi
\ifcase\@tempcnta\tiny\or\scriptsize\or\footnotesize\or\small
 \or\normalsize\or\large\or\Large\or\LARGE\or\huge\else\Huge\fi\fi}
\makeatother

\providecommand{\smaller}[1][1]{\larger[-#1]}

\newcommand{\Smaller}[2][1]{\hbox{\smaller[#1]$#2$}}

\makeatletter \def\@dotsep{.8} \makeatother % почаще точки ...........

\def\red{\color{red}}
\def\blue{\color{blue}}

\newcommand{\Width}[2][\displaystyle]
{\settowidth\Dim{\hbox{$#1#2$}}\texttt{\tiny\blue$($\the\Dim$)$}}
\newcommand{\Height}[2][\displaystyle]
{\settoheight\Dim{\hbox{$#1#2$}}\texttt{\tiny\blue$($\the\Dim$)$}}
\newcommand{\Depth}[2][\displaystyle]
{\settodepth\Dim{\hbox{$#1#2$}}\texttt{\tiny\blue$($\the\Dim$)$}}

\usepackage{hyperref}
\newcommand{\textcenter}[3][Yu.~Brezhnev]{
\pagestyle{myheadings}% задаем текст в колонтитуле
\markboth{\hfill\textsc{#1}\hfill}{\hfill\textsc{#1}\hfill}
\ifdim#2=0mm\textwidth=\textwidth\else\textwidth=#2\fi
\oddsidemargin=\paperwidth%
\addtolength{\oddsidemargin}{-\textwidth}%
\addtolength{\oddsidemargin}{-2in} \oddsidemargin=0.5\oddsidemargin
\evensidemargin=\oddsidemargin
\ifdim#3=0mm\textheight=\textheight\else\textheight=#3\fi
\topmargin=\paperheight%
\addtolength{\topmargin}{-\textheight}%
\addtolength{\topmargin}{-2in}%
\addtolength{\topmargin}{-\headsep}%
\addtolength{\topmargin}{-\headheight}%
\addtolength{\topmargin}{-\footskip}%
\topmargin=0.5\topmargin}

 \textcenter{0mm}{0mm}
\newtheorem{theorem}{Theorem}%[section]
\newtheorem{lemma}[theorem]{Lemma}

\newtheorem{proposition}[theorem]{Proposition}
\newdefinition{definition}{Definition}
\newdefinition{remark}{Remark}
\newdefinition{example}{\sf Example}
\newproof{pf}{Proof}

\allowdisplaybreaks[4]

\def\om{\omega}
\def\T{\hbox{{\footnotesize$\mathrm{T}$}}}
\def\u{{\mathfrak{u}}}
\def\bx{{\boldsymbol{x}}}
\def\bz{{\boldsymbol{z}}}
\def\bs{{\boldsymbol{\mathfrak{s}}}}
\def\bu{{\boldsymbol{\mathfrak{u}}}}
\def\G{{\boldsymbol{\mathfrak{G}}}}
\def\bcQ{{\mathcal{Q}}}
\def\X{\mathrm{X}}
\def\Y{\mathrm{Y}}
\def\Z{\mathrm{Z}}
\def\I{{\s\mathrm{I}}}
\def\II{{\s\mathrm{II}}}
\def\III{{\s\mathrm{III}}}
\def\IV{{\s\mathrm{IV}}}
\def\V{{\s\mathrm{V}}}
\def\VI{{\s\mathrm{VI}}}
\date{today}

\journal{Journal of Differential Equations}

\begin{document}
\begin{frontmatter}

\title{A note on Chudnovsky's Fuchsian equations}

\tnotetext[t1]{Research supported by the Federal Targeted Program under
contract 02.740.11.0238.}

\author{Yurii~V.~Brezhnev}
\ead{brezhnev@mail.ru}

%\address{The Moon}

\begin{abstract}
We  show that  four exceptional Fuchsian equations, each determined by
the four parabolic singularities, known as the Chudnovsky equations,
are transformed into each other by algebraic transformations. We
describe  equivalence of these equations and their counterparts on
tori. The latter are the Fuchsian equations on elliptic curves and
their equivalence is characterized by transcendental  transformations
which are represented explicitly in terms of elliptic and theta
functions.
\end{abstract}


\end{frontmatter}

\tableofcontents \thispagestyle{empty}

\section{Introduction}

\subsection{Chudnovsky equations}
The subject of the present work is the set of four ordinary
differential equations
\begin{alignat}{3}
x\,(x-1)(x+1)\,&\Psi''+{}&(3\,x^2-1)\,&\Psi'+{}&(x+0)\,\Psi=0\,,
\label{1}\\
x\,(x^2+3\,x+3)\,&\Psi''+{}&(3\,x^2+6\,x+3)\,&\Psi'+{}&(x+1)\,\Psi=0\,,
\label{2}\\
x\,(x-1)(x+8)\,&\Psi''+{}&(3\,x^2+14\,x-8)\,&\Psi'+{}&(x+2)\,\Psi=0\,,
\label{3}\\
x\,(x^2+11\,x-1)\,&\Psi''+{}&(3\,x^2+22\,x-1)\,&\Psi'+{}&(x+3)\,\Psi=0
\label{4}\,,
\end{alignat}
reported for the first time by D.~Chudnovsky \& G.~Chudnovsky
\cite{chud1} and considered later more fully in their remarkable work
\cite{chud2}. Once their arising in 1986 it became clear that  list
\eqref{1}--\eqref{4} is quite exceptional and one of the features of
these equations is the fact that these are the only linear ordinary
differential equations (ODEs) of the class
\begin{equation}\label{P}
p\,\Psi''+p'\,\Psi'+(x+A)\,\Psi=0\,,\qquad
p\DEF x\,(x-\alpha)(x-\beta)\,,
\end{equation}
solutions of which are known in terms of known special functions. It is
interesting also to observe that these equations, solvable as they are,
fit no in any currently available algorithmic methods of integration
(over $_2F_1$-extension fields) known in the differential
Picard--Vessiot theory \cite{singer}.

From the Fuchsian standpoint the equations have the parabolic
singularities  at each of the points $x=\{0,\alpha,\beta,\infty\}$,
\ie, Fuchsian exponent differences are equal to zero there. Smirnov, in
his Magister Dissertation \cite{smirnov1} and subsequent work
\cite{smirnov2}, considered equations of the form \eqref{P} and the
question as to their reducibility to a hypergeometric equation by
rational transformations of independent variable $x\mapsto z= R(x)$. He
showed that there are finitely many cases of such reductions and found
one of them. Solutions to equations \eqref{1}--\eqref{4} reduce to the
hypergeometric functions $_2F_1(a,b;c|z)$ indeed. However,
transformations are nontrivial and their complete list was written down
only recently by F.~Beukers \cite{beukers}\footnote{All the cases on
p.~427--428 are correct except for misprint $b(z)^{1/4}\to
b(z)^{\sm1/4}$ and some incorrectness in {\bf case~B} on p.~428. See
also entries (6), (8), and (9) in Tables~12--13 of work \cite{maier}.}.
Arguments of works \cite{beukers} and \cite{chud2} are concerned with
integral recurrences and another (simple) explanation is related to the
fact revealed by A.~Beauville in \cite{beauville}. He found a complete
list of six stable $t$-families of elliptic curves $F_t(\X,\Y,\Z)=0$
over $\mathbb{P}^1(\mathbb{C})$ with only four singular fibres; these
are determined by those $t$-values that degenerate the curve $F_t=0$
into a rational curve (zero genus). In the language of linear ODEs this
entails existence of Heun's equations, all of whose monodromy groups
$\G_t$ are subgroups of the full modular group
$\mathrm{PSL}_2(\mathbb{Z}) \FED\boldsymbol{\Gamma}(1)$  and determine
the zero genus orbifolds $\Hp\!/\G_t$  with four cusps and no elliptic
points. This property was also confirmed by a purely group point of
view in the classification work \cite[Tables 2, 3]{sebbar} and
equivalents of Zagier--Beukers three-term recurrences
\cite[p.~427]{beukers}  were discovered, shortly after Beauville's
list, in Coster's Thesis \cite{coster} as ones associated with
Beauville's curves.

Let us sketch a way of derivation of Beukers' $_2F_1$-reduction
formulae making use of Beauvilles' results. Consider the original
Beauville list \cite[p.~658]{beauville}:
\begin{align}
\X^3+\Y^3+\Z^3+t\,\X\,\Y\,\Z={}&0\,,\tag{I}\label{I}\\
\X\,(\X^2+\Z^2+2\,\Z\,\Y)+t\,\Z\,(\X^2-\Y^2)={}&0\,,\tag{II}
\label{II}\\
\X\,(\X-\Z)(\Y-\Z)+t\,\Z\,\Y\,(\X-\Y)={}&0\,,\tag{III}\label{III}\\
(\X+\Y)(\Y+\Z)(\Z+\X)+t\,\X\,\Y\,\Z={}&0\,,\tag{IV}\label{IV}\\
(\X+\Y)(\X\,\Y-\Z^2)+t\,\X\,\Y\,\Z={}&0\,,\tag{V}\label{V}\\
\X^2\,\Y+\Y^2\,\Z+\Z^2\,\X+t\,\X\,\Y\,\Z={}&0\tag{VI}\label{VI}
\end{align}
and compute  Klein's $J$-invariants for these elliptic curves. We
obtain
\begin{alignat}{5}
&J_{\I} &&=\frac{-1}{12^3}\,\frac{t^3(t-6)^3(t^2+6\,t+36)^3}
{(t+3)^3(t^2-3\,t+9)^3}\,,&\qquad &J_{\IV}
&&=\frac{1}{12^3}\,\frac{(t+2)^3((t+2)^3-24\,t)^3}
{t^3(t+8)(t-1)^2}\,,\notag\\[1ex]
&J_{\II}&&=\frac{4}{27}\,\frac{(t^4-t^2+1)^3} {t^4(t-1)^2(t+1)^2}
\,,&&J_{\V}&&= \frac{1}{12^3}\,\frac{(t^4+16\,t^2+16)^3} {t^2(t^2+16)}
\,,\label{6}\\[1ex]
&J_{\III} &&= \frac{1}{12^3}\,\frac{((t-3)^4-40\,(t^2-3\,t+2))^3}
{t^5(t^2-11\,t-1)} \,,&
&J_{\VI}&&=\frac{-1}{12^3}\,\frac{t^3(t^3+24)^3}
{(t+3)(t^2-3\,t+9)}\,.\notag
\end{alignat}
On the other hand, Klein's $J$ is determined by a classical
hypergeometric Fuchsian equation of the form (H.~Bruns (1875))
\begin{equation}\label{J}
J\,(J-1)\,\Psi_{\!\!\mathit{J\!\!J}}^{}\,+
\frac16(7\,J-4)\,\Psi_{\!J}^{}\,+
\frac{1}{144}\,\Psi=0
\end{equation}
whose monodromy group $\G_{\!J}^{}$ is $\boldsymbol{\Gamma}(1)$. We may
therefore consider formulae \eqref{6} as changes of variables $J\mapsto
t$: $J=R(t)$; each such a change substituted in \eqref{J} must cause
this equation to become the Fuchsian one having monodromy among
Beauville's groups \cite{beauville}, namely, group of a certain
4-punctured sphere. Hence, the resulting ODEs
$\psi''+p(t)\,\psi'+q(t)\,\psi=0$ are solved in terms of
$_2F_1$-solutions to \eqref{J}, that is
\begin{equation}\label{2F1}
\psi(t)=m(t)\cdot{}_2F_1\!\!\!\left(\frac{1}{12},\frac{1}{12};
\frac23\Big|J_k(t)\right),\qquad
k=\mbox{I, II, \ldots, VI},
\end{equation}
where $J_k(t)$ are taken from expressions \eqref{6} and $m(t)$ is an
easily computable multiplier. All this provides a simple way of getting
formulae and results in an equivalent to Beukers' ones
\cite[p.~427--428]{beukers} up to M\"obius transformations of variables
$t$ and renormalization $\Psi\mapsto \psi=m(t)\,\Psi$ which has no
effect on monodromy representations $\G_t$. Doing this, we immediately
reveal (the known fact \cite{beukers}) that $t$-equations for the cases
\eqref{I}, \eqref{VI} coincide and case \eqref{II} is equivalent to
\eqref{V} by a trivial scaling $t\mapsto 4\,\ri\, t$. In the end this
yields the four independent equations which are equivalent to
Chudnovsky's list \eqref{1}--\eqref{4}.

\subsection{Motivation and results}
Transformations between Fuchsian equations of the rational type
$x\mapsto z=R(x)$ are the subject of numerous studies and go beyond
equations with parabolic singularities, hypergeometric reducibility, or
(Heun's) equations with four singular points. Recent results on Heun's
equations have been summarized in work \cite{reiter} (see also
references therein) although first examples appeared already in
\cite{chud2}. However rational transformations are a particular case of
the general algebraic ones which have not yet been considered in the
literature. On the other hand, such a kind substitutions $F(x,z)=0$ may
be thought of as Riemann surfaces and their genera may turn out to be
nontrivial in general. In particular, these considerations allow us to
obtain their parametrizations (uniformizing Hauptmoduln, \ie, principal
moduli in Klein's terminology). These surfaces can split  or not split
to simpler surfaces but genesis and structure of these reducibilities
are presently unknown.

In this note we show that independence of the equations with respect to
Beukers' rational transformations reduces to their common algebraic
equivalence (Sect.~\ref{algchud}); this is done by  certain algebraic
substitutions and leads to very nontrivial results concerning Fuchsian
equations on tori. Correlating such equations with arising algebraic
curves, we obtain  Riemann surfaces that admit the transcendental
representations in form of (mutual) covers of tori (Sect.~\ref{s3}).
Sections~\ref{eqODE} and \ref{mon} contain an additional and more
detailed motivation for algebraic/transcendental equivalence of
equations under study.

\begin{theorem}\label{T1}
Chudnovsky's equations \textup{\eqref{1}--\eqref{4}} and their
counterparts on tori  $($elliptic curves\/$)$ are transformable into
each other by algebraic and  transcendental changes of independent
variables. All the changes are explicitly computable $($listed
below\/$)$ and define the equivalence relations between integrabilities
of these equations.
\end{theorem}

The exclusive character of the list \eqref{1}--\eqref{4} tells us that
these algebraic curves (Table~\ref{capt} and Theorem~\ref{T2}) are also
exclusive since they realize an equivalence of any of Chudnovsky's
equations to any other of them. We also give a treatment to the known
Halphen transformation \cite{poole,halphen} as a transcendental
(bi-single-valued) analog of birational transformations between
polynomial (algebraic) models of an elliptic curve. This allows us to
pass explicitly to associated equations on tori. We tabulate these
equations and their equivalence which is essentially  transcendental
and representable in terms of elliptic functions. This is of special
interest because implicit algebraic dependencies admissible
representations in terms of covers of elliptic tori are very
effectively described through Jacobi's theta-functions. Whilst
Eqs.~\eqref{1}--\eqref{4} define zero genus orbifolds, they explicitly
lead to  Riemann surfaces/orbifolds of higher genera being no
transformations between Chudnovsky's equations. In particular, the
famous Schwarz hyperelliptic curve $y^2=x^8+14\,x^4+1$ appears.

The paper is organized as listed in Contents.

\section{Transformations and equivalences\label{s2}}
In a nutshell,  existence of the above mentioned transitions follows
from the fact that each of groups $\G_t$ in \eqref{6} is a subgroup of
$\boldsymbol{\Gamma}(1)$ and therefore all of these groups are
commensurable each other. Hence it follows that there is a
transformation of algebraic form $F(t_1,t_2)=0$ turning any
$\G_{t_1}$-equation into any other one for $\G_{t_2}$. These algebraic
dependencies are nothing but equalities of $J$-invariants \eqref{6}
between themselves. It turns out that the sought-for algebraic changes
are not always of complicated form coming from a direct equating $J$'s
each other. Below is an example of the most generic case.

\begin{example}\label{E1}
Denote $t$'s for \eqref{III} as $-z$ and $x$ for \eqref{IV} and
consider equality $J_{\III}=J_{\IV}$:
\begin{equation}\label{34}
-\frac{((z+3)^4-40\,(z^2+3\,z+2))^3}
{z^5(z^2+11\,z-1)}
=\frac{(x+2)^3((x+2)^3-24\,x)^3}
{x^3(x+8)(x-1)^2}\,.
\end{equation}
Turning this equation into a polynomial $F(x,z)=0$, we found that it is
irreducible and defines  an algebraic curve of genus $g=5$ .
\end{example}

\subsection{Substitutions} In order to compare Fuchsian equations  it
is convenient to pass to their canonical normal form
$\psi''=\bcQ\,\psi$ because it is unique as against the generic form
$\Psi_{\!\!\mathit{xx}}+p\,\Psi_{\!x}+q\,\Psi=0$. Corresponding linear
transformation $ \Psi\rightleftarrows\psi$ is very well known
\cite{poole,WW,singer} and may be accompanied by a simultaneous  change
of independent variable $x\mapsto z=z(x)$:
\begin{equation}\label{change}
\psi(z)=\sqrt{\frac{dz}{dx}}\,\,\re^{\frac12
\hbox{\tiny{$\ds\int\limits^{\;x}$}}
\!\!p\, dx}_{\mathstrut} \,\Psi(x)\,.
\end{equation}
Then equation for $\psi$ has the form
\begin{equation}\label{star}
\psi_{\mathit{zz}}^{}=\frac12\mbig[8]\{
\frac{z_\mathit{xxx}}{\pow[\textstyle]{z}{x}{3}}-\frac32
\frac{\pow{z}{\mathit{xx}}{2}}{\pow[\textstyle]{z}{x}{4}}
+\frac{1}{\pow[\textstyle]{z}{x}{2}}
\Big( p_x^{}+\Mfrac12\,p^2-2\,q\Big)\!
\mbig[8]\}\psi\,.
\end{equation}
Intermediate transformations $x\mapsto x'\mapsto x''\mapsto\cdots
\mapsto z$  are allowable but the number of such changes and their
orders, including inverse transformations, are immaterial for ultimate
answer $x\rightarrow z$; this formula has an invariant
characterization.

In practice, when the change $x\mapsto z$ has been given in form of
implicit equation $F(z,x)=0$,  it is useful to have an effective
formula for transition to the normal form $\psi''=\bcQ(z)\,\psi$, where
primes, as always in the sequel, signify the derivatives with respect
to independent variable entering into coefficient of the proper
$\psi$-equation. As usual, when transforming linear ODEs the Schwarz
derivative does constantly appear  and we use the standard notation for
this object:
$$
\{f,z\}\DEF
\frac{f_{\!\mathit{zzz}}}{f_z}-
\frac32\frac{\pow{f}{\mathit{zz}}{2}}{\pow{f}{z}{2}}\,.
$$
With use of this notation we can rewrite the transformation above in
form of the following computational rule.

\begin{lemma}\label{L1}
Let coefficients of equation
\begin{equation}\label{pq}
\Psi_{\!\!\mathit{xx}}+p\,\Psi_{\!x}+q\,\Psi=0
\end{equation}
be arbitrary $($rational, algebraic, or transcendental\/$)$
differentiable functions of $x$. Then linear change \eqref{change} and
the change of variables $x\mapsto z$ defined by the rule $F(z,x)=0$
transform Eq.~\eqref{pq} into the following canonical form:
\begin{equation}\label{z}
\psi''=\frac12\,\bcQ(z)\,\psi\,,
\end{equation}
$$
\bcQ(z)=
\frac{\pow{F}{z}{2}}{\pow{F}{x}{2}}
\mbig[7](p_x^{}+\frac12\,p^2-2\,q+\{F,x\}\mbig[7])
\!-\frac{F_z}{F_x}\,p_z^{}
-\{F,z\}+
3\,\frac{F_z}{F_x}\!
\left(\ln\!\frac{F_z}{F_x}\right)_{\!\!\!\!\mathit{xz}}\,,
$$
where objects   $\{F,x\}$, $\{F,z\}$ are understood as the partial
Schwarz derivatives and expression for $\bcQ(z)$ should be computed
modulo $\big\langle F(z,x)\big\rangle$.
\end{lemma}

\begin{pf} Compute the derivatives $z_x$, $z_\mathit{xx}$, and
$z_\mathit{xxx}$ appearing in \eqref{star} according to the rules like
$$
z_x=-\frac{F_x}{F_z}\,,\qquad
z_\mathit{xx}=-\left(\frac{F_x}{F_z}\right)_x-
\left(\frac{F_x}{F_z}\right)_zz_x=-\frac{F_\mathit{xx}}{F_z}
+2\frac{F_{\mathit{xz}}F_x}{\pow{F}{z}{2}}-
\frac{F_{\mathit{zz}}\pow{F}{x}{2}}{\pow Fz3}\,,\qquad\ldots\,.
$$
Express third derivatives $F_{\mathit{xxx}}$ and $F_{\mathit{zzz}}$ via
partial Schwarzians $\{F,x\}$,  $\{F,z\}$. Taking into account that $p$
may be an algebraic function $p(x,z)$, we replace the complete
derivative $p_x$ presented in \eqref{star} with the following object:
$$
p_x\mapsto p_x-\frac{F_x}{F_z}\,p_z\,.
$$
Simplifying the result, one arrives at  formula for $\bcQ(z)$ above.
\end{pf}

The rule \eqref{z} is convenient to use because its last term vanishes
if the dependence  $F(z,x)=0$ has a split  form $X(x)=Z(z)$, which is
frequently our case. Such form simplifies calculations of genera of
curves and  reduces considerably computation tasks when the polynomial
operation modulo $\big\langle F(z,x)\big\rangle$ has been applied to
the answer $\bcQ(z)$. We shall exploit this lemma throughout the work.

\subsection{On equivalence  of 2nd order linear ODEs\label{eqODE}}
The main motivation for study of transformations between equations
under considerations is the fact that the simple or complicated
Fuchsian (not necessarily) equations may be transformed into very
simple equations with avail of far non-obvious
rational/algebraic/transcendental substitutions.

\begin{proposition}
Any two linear 2nd order ODEs
\begin{equation}\label{12}
\Psi_{\!\!\mathit{xx}}=\bcQ(x)\,\Psi\,,\qquad
\psi_{\mathit{zz}}^{}=\widetilde\bcQ(z)\,\psi
\end{equation}
can be transformed into each other by a point transformation
$z=\Xi(x)$.
\end{proposition}
\begin{pf}
Linearity and normality of both of Eqs.~\eqref{12} implies the linear
relation between $\Psi$ and $\psi$, \eg, $\psi=m\,\Psi$, with
$m=\sqrt{z_x}$, where dependence $z=\Xi(x)$ is as yet unknown. Hence
\begin{equation}\label{temp}
\frac{dx}{\Psi^2}=\frac{dz}{\psi^2}\,.
\end{equation}
Whatever the solution $\Psi=\Psi(x)$ is chosen, we can construct the
second linearly  independent one by Liouville's formula
$$
\Psi\int\!\!\frac{dx}{\Psi^2}\,.
$$
Therefore $\int\Psi^{\sm2}\,dx$ is always a certain ratio of two linear
independent solutions to the $\Psi$-equation and this ratio will be the
same for the $\psi$-equation; the ratio depends only on point $x$. We
thus have, instead of \eqref{temp},
$$
\frac{\Psi_2(x)}{\Psi_1(x)}=\frac{\psi_2^{}(z)}{\psi_1^{}(z)}
$$
and this relation constitutes an implicit form of the sought-for
dependence $z=\Xi(x)$.
\end{pf}

As can well be imagined, such a construction is useless in general
because it requires the knowledge of  integrals. Since the
$\{\psi,\,\Psi\}$'s are chosen to be arbitrary the general equivalence
of Eqs.~\eqref{12} can be rewritten in  form of the following bilinear
relation:
\begin{equation}\label{abcd}
x\rightleftarrows z:\qquad \mathrm{A}\,\Psi_1(x)\,\psi_1(z)+
\mathrm{B}\,\Psi_1(x)\,\psi_2(z)+
\mathrm{C}\,\Psi_2(x)\,\psi_1(z)+\mathrm{D}\,\Psi_2(x)\,\psi_2(z)=0
\end{equation}
with free constants $(\mathrm{A}:\mathrm{B}:\mathrm{C}:\mathrm{D})$.

In the majority of cases integrals of linear ODEs belong to
differential fields which are different from those to which the
coefficients $\bcQ(x), \tilde\bcQ(z)$ belong. It is not at all obvious
a priori, then, that Chudnovsky's equations admit situations when the
family \eqref{abcd} has algebraic representatives $F(x,z)=0$, whereas
all the functions $\Psi_1$, $\Psi_2$, $\psi_1^{}$, and $\psi_2^{}$ are
expressed solely in terms of non-algebraic hypergeometric
${}_2F_1$-transcendents (Beukers' list \cite{beukers}).

\begin{definition}
We shall call linear ODEs \eqref{12} algebraically equivalent if they
are transformable into each other by some algebraic dependence
$F(x,z)=0$.
\end{definition}

\begin{remark}
It is not difficult to see that algebraic equivalence defines an
equivalence relation since it satisfies the symmetry, reflection, and
transitivity  properties. We do not use the separate term for rational
equivalence, \eg, $z=R(x)$, because inversions of the rational function
$R(x)$ and the change $\psi=\sqrt{R'(x)}\,\Psi$ always lead to
algebraic functions. The transcendental equivalence is always
available; this is formula \eqref{abcd}. However in Sect.~\ref{trans}
we shall exhibit examples---Chudnovsky's equations on tori---when
equivalence is transcendental but it is simpler than the most general
one defined by this formula. It may be also mentioned here that
algebraic equivalence is a simplest but nontrivial kind of
equivalences.
\end{remark}

\subsection{Remarks on monodromy groups\label{mon}}
Yet another point that should be mentioned is the fact that
Eqs.~\eqref{1}--\eqref{4} provide the next nontrivial (after a
hypergeometric equation) examples of what is called presently the
monodromy groups of finite genus. Recall that this property implies
that  function  $x=\chi(\tau)$ defined by inversion of the ratio
\begin{equation}\label{ratio}
\tau=\frac{\Psi_1(x)}{\Psi_2(x)}
\end{equation}
is a single-valued analytic function of variable $\tau$ everywhere in
the domain of its existence on the plane $(\tau)$. As usual,  the
closed paths\footnote{If coefficient is an algebraic function
$\bcQ(x,y)$ belonging to irrationality $F(x,y)=0$ then the closure of a
path is defined by the value of the pair $(x,y)$ coinciding with the
initial one $(x_{\s0},y_{\s0})$ \cite{poincare}.} on the plain $(x)$
entail transformations
$$
\bigg(\begin{matrix} \Psi_1\\\Psi_2\end{matrix}\bigg)\mapsto
\bigg(\begin{matrix} a&b\\c&d\end{matrix}\bigg)
\bigg(\begin{matrix} \Psi_1\\\Psi_2\end{matrix}\bigg)
$$
which form a finitely generated group $\G$ (the monodromy group)
\cite{poincare,yoshida} as equation is of Fuchsian class. Then the
$\tau$-plane is covered by domains containing $\G$-nonequivalent points
and pairwise equivalent points on boundaries of the domains. If these
domains form a set of non-overlapping circle polygons with finitely
many number of sides each (Poincar\'e polygons) then identifications of
these sides determine the standard topological characteristics of the
polygon---the genus \cite{ford}; in doing so, the function $\chi(\tau)$
becomes single-valued by construction. For brevity, we shall use
terminological shorthand the monodromy and genus of the monodromy as
synonyms to the monodromy group and genus of the Poincar\'e polygon
representing the group. Being a matrix group from
$\mathrm{SL}_2(\mathbb{C})$, it has an exact representation by an
automorphism group of the (automorphic) function $\chi$ and hence we
re-denote this group as $\G_x$:
$$
\chi\Big(\Mfrac{a\,\tau+b}{c\,\tau+d}\Big)=\chi(\tau)\quad
\hence\quad \mathrm{Aut}\,\chi(\tau)\FED\G_x\,.
$$

If $\bcQ$ is a rational function of $x$ then the monodromy has a zero
genus \cite{ford,poincare}. If $\bcQ=\bcQ(x,y)$ is an algebraic
function belonging to irrationality $F(x,y)=0$ then the genus, by
construction, coincides with topological genus of this curve.  We shall
also meet Fuchsian equations wherein $\bcQ$  is an elliptic
(transcendental) function $\bcQ(\u)$. In this case, genus of monodromy
is, again by construction, equal to unity.

Explicit $\chi(\tau)$-expressions  for equations under consideration
can be found in works \cite{sebbar,maier} and Eqs.~\eqref{1}--\eqref{3}
are related to the classical modular equations as particular cases; the
most exhaustive literature and systematic lists of results concerning
this subject can be found in \cite{maier}.

The arbitrary substitutions $x\mapsto z$ destroy in general the
property of monodromies to have finite genus but it is clear that any
single-valued rational/transcendental change $z=R(x)$ will
automatically yield equation \eqref{z} with the monodromy $\G_z$ known
to be Fuchsian, \ie, of finite genus, if the monodromy $\G_x$ was of
the same kind. However, this is somewhat trivial way to construct new
interesting equations because they will have in general the complicated
algebraic coefficients $\bcQ(x,z)$. As we shall see, the theory of
Chudnovsky's equations provides a large number of nontrivial situations
when rational functions $\bcQ(x)$ with  zero genus monodromies go into
rational functions $\widetilde{\bcQ}(z)$ again, whereas the
substitutions themselves have nontrivial genera. Similarly,  the unity
genera pass to the unity ones (punctured tori; sects.~\ref{51} and
\ref{52}). In other words, genus of manifold on which ODE has been
defined, genus of its monodromy group, and genus of the substitution
are \textit{not one and the same}. Therein lies an essential feature of
algebraic equivalence of Chudnovsky's equations and Fuchsian
monodromies at all.

\subsection{Genera of substitutions\label{genera}}

Turning back to  Eqs.~\eqref{1}--\eqref{4}, let us tabulate their
canonical forms for further reference. We apply the `linear part' of
Lemma~\ref{L1} (\ie, independent variable is not changed) and derive
that normal forms to equations \eqref{1}--\eqref{3} become respectively
\begin{alignat}{2}
\psi''&=-\frac14\,\frac{(x^2+1)^2}{x^2(x-1)^2(x+1)^2}\,
\psi\,,\tag{$1'$}\label{1'}\\[1ex]
\psi''&=-\frac14\,\frac{(x+1)(x+3)(x^2+3)}{x^2(x^2+3\,x+3)^2}\,
\psi\,,\tag{$2'$}\label{2'}\\
\intertext{and}
\psi''&=-\frac14\,\frac{x^4+8\,x^3+72\,x^2-64\,x+64}
{x^2(x-1)^2(x+8)^2}\,\psi\,. \tag{$3'$}\label{3'}\\
\intertext{The normal form for Eq.~\eqref{4} can be obtained
analogously, however, by way of illustration of the two last sentences
in the previous section, we apply Lemma~\ref{L1} in its full generality
and obtain that the change \eqref{34} transforms Eq.~\eqref{3'}
(replacing $z$ with $x$ again) into the following equation}
\psi''&=-\frac14\,\frac{x^4+12\,x^3+134\,x^2-12\,x+1}
{x^2(x^2+11\,x-1)^2}\,\psi\,.\tag{$4'$}\label{4'}
\end{alignat}
This is exactly the normal form to  equation \eqref{4}.  We shall refer
to Eqs.~\eqref{1'}--\eqref{4'} as Chudnovsky's equations as well. To
avoid confusion, we also adjust Beauville's $t$-parameters in \eqref{6}
in order to make exact correlation of these $J$-invariants with list
\eqref{1'}--\eqref{4'}:
\begin{alignat}{5}
\eqref{I}:&\quad t=-3\,(x+1)\,,\qquad&
\eqref{II}:&\quad t=x\,,\qquad&
\eqref{III}:&\quad t=-x\,,\notag\\
\eqref{IV}:&\quad t=x\,,\qquad&
\eqref{V}:&\quad t=4\,\ri\,x\,,\qquad&
\eqref{VI}:&\quad t=-3\,(x+1)\,.\notag
\end{alignat}
Invariants \eqref{6} then read
\begin{alignat}{5}
&J_{\I}&&=\frac{1}{64}\,\frac{(x+1)^3(x+3)^3(x^2+3)^3}
{x^3(x^2+3\,x+3)^3}\,,\qquad& &J_{\IV}
&&=\frac{1}{12^3}\,\frac{(x+2)^3((x+2)^3-24\,x)^3}
{x^3(x+8)(x-1)^2}\,,\notag\\[1ex]
&J_{\II}&&=\frac{4}{27}\,\frac{(x^4-x^2+1)^3} {x^4(x-1)^2(x+1)^2}
\,,&\qquad &J_{\V}&&= \frac{1}{108}\,\frac{(16\,x^4-16\,x^2+1)^3}
{x^2(x-1)(x+1)}
\,,\label{JJ}\\[1ex]
&J_{\III} &&=\frac{-1}{12^3}\,\frac{((x+3)^4-40\,(x^2+3\,x+2))^3}
{x^5(x^2+11\,x-1)} \,,&\qquad
&J_{\VI}&&=\frac{1}{64}\,\frac{(x+1)^3(9\,(x+1)^3-8)^3}
{x\,(x^2+3\,x+3)}\notag
\end{alignat}
and generate, by Lemma~\ref{L1} applied to the `hypergeometry'
\eqref{J}, the list \eqref{1'}--\eqref{4'} as follows:
$$
\eqref{1'}\;\scalebox{-1}[1]{$\mapsto$}\;\{\mathrm{\ref{II}},
\mathrm{\ref{V}}\}\,,\qquad
\eqref{2'}\;\scalebox{-1}[1]{$\mapsto$}\;\{\mathrm{\ref{I}},
\mathrm{\ref{VI}}\}\,,\qquad
\eqref{3'}\;\scalebox{-1}[1]{$\mapsto$}\;\mathrm{\ref{IV}}\,,\qquad
\eqref{4'}\;\scalebox{-1}[1]{$\mapsto$}\;\mathrm{\ref{III}}\,.
$$

\begin{theorem}\label{T0}
Let algebraic equivalence be generated by identifying Klein's
$J$-invariants \eqref{JJ}. Then Chudnovsky's equations
\textup{\eqref{1'}--\eqref{4'}} are algebraically equivalent with
respect to substitutions whose genera are presented in
Table~$\ref{capt}$.
\begin{table}
\setlength{\doublerulesep}{0.5pt}
\renewcommand{\arraystretch}{1.1}
\begin{center}\begin{tabular}{|r||c|c|c|c|c|c|} \hline\label{tab}
&\ref{I}$\eqref{2'}_x$&\ref{II}$\eqref{1'}_x$&
\ref{III}$\eqref{4'}_x$&\ref{IV}$\eqref{3'}_x$&
\ref{V}$\eqref{1'}_x$&\ref{VI}$\eqref{2'}_x$\\
\hline\hline
\ref{I}$\eqref{2'}_z$&$0_{(12)}$&5&5&$0_{(4)}$&5& $\{1, 1, 1, 1\}_
\mathrm{\s E}^{}$\\
\hline \ref{II}$\eqref{1'}_z$&&$0_{(5)}$&5&$\{1, 1, 1\}_ \mathrm{\s
C}^{}$&$0_{(2)}$, $\{1, 1\}_
\mathrm{\s L}^{}$&5\\
\hline \ref{III}$\eqref{4'}_z$&&&$0_{(4)}$&5&5&5\\ \hline
\ref{IV}$\eqref{3'}_z$&&&&$0_{(3)}$, $1_\mathrm{\s E}$&$1_ \mathrm{\s
C}$, 3&0, 4\\ \hline \ref{V}$\eqref{1'}_z$&&&&&$0_{(3)}$, 3&5\\ \hline
\ref{VI}$\eqref{2'}_z$&&&&&&$0_{(3)}$, 4\\ \hline
\end{tabular}\end{center}
\vspace{-.5em} \caption{\label{capt}
Genera of curves realizing
algebraic equivalencies of Chudnovsky's equations.}
\end{table}
\end{theorem}

\begin{pf}
Consider equalities of $J$-invariants \eqref{JJ}: $J_k(x)=J_n(z)$ and
take, \eg, the case $J_{\II}(z)=J_{\III}(x)$. It corresponds to a table
record on intersection of line \mbox{\ref{II}$\eqref{1'}_z$} and column
\mbox{\ref{III}$\eqref{4'}_x$}. Converting this equation into a
polynomial $F(x,z)=0$, we establish that it is not reducible over
$\mathbb{C}$. Since this polynomial represents the equality of one and
the same quantity---Klein's invariant $J$---it ensures the mutual
equivalence of Chudnovsky's
Eqs.~$\eqref{1'}_z\rightleftarrows\eqref{4'}_x$; of course, this can be
checked by a straightforward  application of Lemma~\ref{L1}.
Computation of genus $g$ by the Riemann--Hurwitz formula gives $g=5$.
Such an irreducibility is not a common rule and we take, as a second
instance, equation $J_{\II}(z)=J_{\V}(x)$. Corresponding polynomial
$F(x,z)=0$ splits into several  components
$$
\begin{aligned}
F(x,z)&=\big((z-1)^2+4\,x^2z\big)\big( (z+1)^2-4\,x^2z\big)\\
&\quad\;\times\!\!\!\big(16\,(x^4-x^2)(z^4-z^2)-1\big)
\big(16\,(x^4-x^2)(z^2-1)+
z^4\big)=0
\end{aligned}
$$
and direct computations (Lemma~\ref{L1}) show that each of them does
realize an algebraic equivalence of Eq.~\eqref{1'} with itself. The
first two components are the rational algebraic curves; their genera
are equal to zero. This point has been designated in the table as
$0_{(2)}$; subscript stands for a number of rational curves and trivial
substitution $x=z$ is taken into account for diagonal cases. The two
unities $\{1,1\}_\mathrm{\s L}^{}$ in the entry mean that the two
remaining polynomials determine curves of genus $g=1$ and each of the
curves is isomorphic to a lemniscate (\textsc{l}); \ie, their  Klein's
$J$-invariants are equal to 1. The symbol $1_\mathrm{\s E}$ designates
a curve isomorphic to the equi-anharmonic (\textsc{e}) curve ($J=0)$
and $1_\mathrm{\s C}$ (Chudnovsky) does the curve with invariant
$J=\frac{13^3}{2^23^5}$. Other entries of the table are processed in a
similar manner and all the curves are defined over
$\mathbb{Q}(\sqrt{5})$ or $\mathbb{Q}(\ri\,\sqrt{3})$ at most
(splitting fields of polynomials $x^2+11\,x-1$ and $x^2+3\,x+3$). In
all the cases curves of the same genus differ from one another and can
be rather complicated. The empty boxes are filled by symmetry.
\end{pf}

\begin{remark}\label{R0}
We do not have an explanation of   `unpredictable' distribution of
genera in  Table~\ref{capt} or explanation as to why \textit{each} of
irreducible components does indeed represent an algebraic equivalence
of the list \eqref{1'}--\eqref{4'}. This is, perhaps, a quite
nontrivial task because it touches  upon the problem of construction of
all the algebraic equivalences. This goes far beyond the scope of the
present work and, as  mentioned in Introduction, no the general theory
of algebraic transformations has yet  been developed.
\end{remark}

\section{Algebraic equivalence of Chudnovsky equations\label{algchud}}

\subsection{Automorphisms and their consequences\label{aut}}
Table~\ref{capt} shows that there are transformations of the same
Chudnovsky equation into itself and these are defined not only through
the trivial change $z=x$. Non-obvious examples appear even in the class
of linear fractional  substitutions. For example, the zero genus family
of automorphisms \eqref{2'}$~\rightleftarrows$~\eqref{2'} contains the
transformation
$$
(\varepsilon-1)(x+z) = x\,z+3\,,\qquad
\varepsilon\DEF \re^{\frac23\pi\ri}_{\mathstrut}\,.
$$

The mere fact that such transformations do exist is not surprising.
Well-known examples are the modular Jacobi--Schl\ae{}fli--Sohnke
relations between Legendre's moduli $x=k^2(\tau)$ and $z=k^2(N\tau)$.
Diagonal cases $\eqref{1'}~\rightleftarrows~\eqref{1'}$ are thus
particular analogs of this classical family and other diagonal entries
provide certain modular equations belonging to their Beauville
monodromy groups. For example, the right lower box of the table
contains a genus $g=4$ modular equation for Beauville's VI-group
\cite[p.~658]{beauville}, which is conjugate to group
$\boldsymbol{\Gamma}(3)$ \cite{sebbar,maier,coster}. It is of more
interest that transformations of such a kind lead to other interesting
consequences. They are concerned with rational and elliptic
automorphisms. We present here consequences of only two illustrative
examples.

\begin{example}\label{E2}
Let us consider one of the zero genus diagonal quadratic automorphisms
\begin{equation*}\label{ex}
\eqref{3'}_x\rightleftarrows\eqref{3'}_z:\qquad z\,x\,(z+x+6)=8\,.
\end{equation*}
We can parametrize this dependence by rational functions
\begin{equation*}\label{t1}
x=\frac{(\T-1)^2}{\T+1}\,,\qquad
z=\frac{8}{\T^2-1}
\end{equation*}
and may consider $x=x(\T)$ as a change of variable $x\mapsto \T$ in
Chudnovsky's equation \eqref{3'}. Insomuch as rational uniformizer $\T$
itself is always  a rational function of coordinates $(x,z)$, that is
$\T=R(x,z)$, and Fuchsian equation \eqref{3'} has a correct accessory
parameter \cite{chud2}, the transformed $\T$-equation will be of the
same property. Of course, both of these substitutions will yield the
same $\T$-equation. It turns out that equations generated by this way
become new Fuchsian ones and renormalization of  $\T$ can impart them
better (canonical) form. We therefore replace the last parametrization
with this one:
$$
x=2\,\frac{(\T-1)^3}{1-\T^3}\,,\qquad
z=6\,\frac{\T+\varepsilon\,\T+\varepsilon}{1-\T^3}\,\T-2
$$
and derive that $\T$-equation has a very elegant form indeed:
$$
\psi''=\frac{-9\,\T^4}{(\T^6-1)^2}\,\psi\,.
$$
It is an equation of the same kind as Chudnovsky's ones, with the
difference that it has six parabolic singularities at points
$\T=\pm\{1,\varepsilon,\varepsilon+1\}$. One can show, with use of some
manipulations by Jacobi's $\vartheta$-constant series (this is not a
subject matter of the present work), that the $\tau$-representation for
the corresponding Hauptmodul $\T=\T(\tau)$ has the form
$$
\T(\tau)=\varepsilon\,
\frac{\ri\,\vartheta_3^2(\tau)+\sqrt{3}\,\vartheta_3^2(3\tau)}
{\ri\,\vartheta_3^2(\tau)-\sqrt{3}\,\vartheta_3^2(3\tau)}\,,
$$
where the standard Jacobi theta-constant $\vartheta_3(\tau)$ is defined
by the series \cite{weber}
$$
\vartheta_3(\tau)\DEF\sideset{}{_k}\sum_{-\infty}^{+\infty}\!
\re^{k^2\pi\ri\tau}\,.
$$
As it follows from the last $\psi$-equation this $\T(\tau)$ satisfies
the equation
$$
\frac{\{\T,\tau\}}{\dot\T{}^2}=\frac{-18\,\T^4}
{(\T^6-1)^2}\,.
$$
\end{example}

\begin{example}[Non-diagonal automorphisms]\label{E3}
Such automorphisms are possible only for equations \eqref{1'} and
\eqref{2'}. Table~\ref{capt} tells us that these curves are rational or
elliptic ones and  they are isomorphic as curves under coinciding
genera. Choosing the simplest representatives, we obtain in these
cases:
\begin{alignat}{2}
\eqref{2'}_x\rightleftarrows\eqref{2'}_z:\qquad&
\big\{(x+1)^3-1\big\}\big\{(z+1)^3-1\big\}=1\,,&&\quad
g=1 \quad(J=0)\notag\\
\intertext{and}
\eqref{1'}_x\rightleftarrows\eqref{1'}_z:\qquad&
16\,(x^4-x^2)\,(z^4-z^2)=1\,,&&\quad
g=1 \quad(J=1)\,,\label{11}\\
&4\,x^2\,z=(z+1)^2\,,&&\quad g=0\,.\notag
\end{alignat}
The latter case shows that variable $z$ is a perfect square and we can
put  $z=\T^2$, where $\T$ is a uniformizer for this rational curve.
Computing Fuchsian $\T$-equation, we get
$$
\psi''=-\frac14\,\frac{\T^8+14\,\T^4+1}{\T^2(\T^4-1)^2}\,\psi\,,
$$
that is yet another (known \cite{maier}) Fuchsian ODE with six
parabolic singularities.
\end{example}

Eq.~\eqref{11} will be used in Sect.~\ref{trans} when considering
equations on tori.

\begin{remark}\label{R3}
Notwithstanding the fact that one has six Beauville's curves and just
four Chudnovsky's equations we cannot discard some two of curves
\eqref{I}, \eqref{II}, \eqref{V}, \eqref{VI} as excessive. For example,
if we cut out the left upper $(4\!\times\!4)$-box from
Table~\ref{capt}, we would lose many transformations: all the non-zero
genus automorphisms $\eqref{1'}\rightleftarrows\eqref{1'}$,
$\eqref{2'}\rightleftarrows\eqref{2'}$, genus 3 transformation
$\eqref{3'}\rightleftarrows\eqref{1'}$, etc. In this respect
Beauville's list is independent of Chudnovsky's one.
\end{remark}

\begin{remark}\label{R01}
Automorphisms coming from Table~\ref{capt} are not the only possible
ones. This results from the fact that some of Chudnovsky's Hauptmoduln
are expressed via the classical Jacobi theta-constants
\cite{maier,sebbar} which are algebraically related to Legendre's
modulus $k^2(\tau)$. The latter, as is well known, has a lot (infinite)
of algebraic dependencies with itself $k^2(q\,\tau)$ under
$q\in\mathbb{Q}$. Thus, all the possible automorphisms will be analogs
of the classical modular families mentioned above.
\end{remark}

\subsection{Equivalences and integrability of the list
\protect\eqref{1}--\protect\eqref{4}} Although automorphisms produce a
large number of nontrivial curves,  we shall consider further
equivalences of only pairwise distinct  equations.

\begin{theorem}\label{T2}
Algebraic equivalence defined by Table~$\ref{capt}$ is closed if
automorphisms are excluded form consideration. This closedness is
determined by the family of genus $5$ irreducible algebraic curves:
\begin{alignat}{10}
\eqref{1'}_z\rightleftarrows\eqref{2'}_x:&& J_{\II}(z)&=J_{\I}(x)\,,&
J_{\V}(z)&=J_{\I}(x)\,,\notag\\[1ex]
&&J_{\II}(z)&=J_{\VI}(x)\,,&
J_{\V}(z)&=J_{\VI}(x)\,,\notag\\[1ex]
\eqref{1'}_z\rightleftarrows\eqref{4'}_x:&& J_{\II}(z)&=J_{\III}(x)\,,&
J_{\V}(z)&=J_{\III}(x)\,,\notag\\[1ex]
\eqref{2'}_z\rightleftarrows\eqref{4'}_x:&& J_{\I}(z)&=J_{\III}(x)\,,&
J_{\VI}(z)&=J_{\III}(x)\,,\notag\\[1ex]
\eqref{3'}_z\rightleftarrows\eqref{4'}_x:&&
J_{\IV}(z)&=J_{\III}(x)\,,\notag \intertext{and  two exceptional cases
determined by the $($canonical\/$)$ representatives of minimal genera.
These are the elliptic curve}
\eqref{3'}_z\rightleftarrows\eqref{1'}_x:&& \quad
x^2-x^4&=\frac{16\,(z-1)}{z^3(z+8)}&& \hspace{-2em}
\mbig[8](J=\frac{13^3}{2^2\,3^5}\mbig[8])\label{J123}\\
\intertext{and the zero genus one}
\eqref{3'}_z\rightleftarrows\eqref{2'}_x:&&\quad
(x+1)^3&=\frac{(z+2)^3}{(z+8)(z-1)^2} \,.\label{zero}
\end{alignat}
\end{theorem}
\begin{pf}
Except for equalities of Klein's $J$-invariants we should check the
transitivity of relations under consideration. All the irreducible
relations $J_k(x)=J_n(z)$ are listed in the first family of the curves
above. Reducible cases, according to Table~\ref{capt}, are
\eqref{3'}~$\rightleftarrows$\eqref{1'},
\eqref{3'}~$\rightleftarrows$~\eqref{2'}, and each of automorphisms
$J_k(x)=J_k(z)$. Consider, \eg, transitivity
$$
\eqref{2'}_x\rightarrow\eqref{3'}_\bz\,,\qquad
\eqref{3'}_\bz\rightarrow\eqref{4'}_z\,.
$$
For each of the \eqref{2'}~$\rightarrow$~\eqref{3'}-relations we expect
to get the $\eqref{2'}_x\rightarrow\eqref{4'}_z$-relation coinciding
with one of the two table curves of genus 5. There are two sets of the
\eqref{2'}~$\rightarrow$~\eqref{3'}-transformations: $\{0\}_4^{}$ and
$\{0,4\}$. Take a curve from the first set, \eg, the curve
\eqref{zero}:
$$
(x+1)^3=\frac{(\bz+2)^3}{(\bz+8)(\bz-1)^2}
$$
and supplement it with the unique
\eqref{3'}~$\rightarrow$~\eqref{4'}-curve (see Example~\ref{E1})
$$
\frac{(\bz+2)^3((\bz+2)^3-24\,\bz)^3}
{\bz^3(\bz+8)(\bz-1)^2}=-\frac{((z+3)^4-40\,(z^2+3\,z+2))^3}
{z^5(z^2+11\,z-1)}\,.
$$
Elimination of variable $\bz$ from these two equations produces not a
new relation but irreducible curve $J_\I(x)=J_\III(z)$; more precisely,
cube of the curve $\eqref{2'}\rightleftarrows\eqref{4'}$. Other
elements of the sets and transitivity
$\eqref{1'}_x\rightarrow\eqref{3'}_\bz\rightarrow\eqref{4'}_z$ are
checked in a similar manner. The simplest representative of the
\eqref{1'}~$\rightleftarrows$~\eqref{3'}-equivalence with a minimal
genus is the curve \eqref{J123}. By virtue of irreducibility and
transitivity we may also leave single representatives for each of the
cases in the $g=5$ family above\footnote{This rises the question as to
a correlation between these curves. In particular, whether they are
isomorphic or not? }. As for automorphisms, many of them, rather
`exotic' as they are, preserve the closedness of Table~\ref{capt}.
However, there are exceptions. It will suffice to point out at least
one counterexample. This is equivalence $J_\III\to J_\I$ followed by
application of a $g=4$ automorphism coming from the
$(J_\VI\rightleftarrows J_\VI)$-curve. Corresponding transformation
$\eqref{4'}\to\eqref{2'}\to\eqref{2'}$ leads to a cumbersome curve
\hbox{$F(\stackrel{36}{x},\stackrel{36}{z})=0$} of genus $g=25$
(computation is very nontrivial; 36 is a degree of the curve in both
the variables).
\end{pf}

We observe in passing that case $\eqref{4'}$ holds an exceptional
position among other equations $\eqref{1'}$--$\eqref{4'}$ since it is
transformed into other ones only by means of the most complicated
changes. The zero genus automorphisms of this equation, apart from
trivial ones $x=z$ and $x\,z+1=0$, are rather non-obvious and
cumbersome (not displayed here). Hauptmodul for this equation is also
very non-standard \cite[Table~3]{sebbar}. We also see that genus 5
transformation $\eqref{1'}\rightleftarrows\eqref{2'}$  can be
represented as a composition of the simple rational \eqref{zero} and
elliptic curve  \eqref{J123}.

Let us consider the question on integrability of Chudnovsky's
equations. The hypergeometric series converges only in a unite circle,
which is why it would be more convenient to have solutions expressed
not in terms of Beukers' ${}_2F_1$-list but in terms of special
functions associated with the hypergeometric equation. These are
Legendre's complete elliptic integrals $\ellK(k)$, $\ellK'(k)$
\cite{WW,poole} or general Legendre's $P,Q$-functions solving the
equation \cite[Sect.~15$\boldsymbol{\cdot}$5]{WW}
\begin{equation}\label{leg}
(1-s^2)\,Y_{\mathit{ss}}-2\,s\,Y_s+
\big\{\nu(\nu+1)-\mu^2(1-s^2)^{\sm1}\big\}Y=0\,.
\end{equation}

Integrability of Eqs.~\eqref{1}--\eqref{4} in terms of the integrals
above is obvious because the first Chudnovsky's equation \eqref{1} is
in effect equation for  a square root of the standard Legendre's
elliptic modulus $k^2(\tau)=x^2$ defined by the classical equation
\cite{WW}
$$
\eqref{1}\quad\Leftrightarrow\quad\frac{d}{dk}\!\!\left( k\,(1-k^2)\,
\frac{d\,\psi}{dk}\right)=k\,\psi\,,\qquad
\psi=\big\{\ellK(x),\,\ellK'(x)\big\}\,.
$$
It is common knowledge that there are cases when the $_2F_1$-series
admits the quadratic rational transformations and the generic
hypergeometric equation then reduces to the two-parametric equation
\eqref{leg} \cite[Sect.~3.1--2]{bateman2}. This is indeed the case for
equations under question and some simple arguments show that one of the
reductions is  $(\nu,\mu)=({-}\frac13, 0)$.

\begin{proposition}\label{C2}
All the equations \textup{\eqref{1}--\eqref{4}} are integrable in terms
of Legendre's integrals $K$, $K'$ or functions $P_{\!\sm1\!/\!3}$,
$Q_{\!\sm1\!/\!3}$.
\end{proposition}
\begin{pf}
It will suffice to integrate one equation of the list
\eqref{1'}--\eqref{4'}. We take Eq.~\eqref{2'} and derive that it is
transformed into Eq.~\eqref{leg} and the `hypergeometry' \eqref{J} as
follows
$$
J=\frac14\,\frac{(4\,s-5)^3}{(s^2-1)(s+1)}\qquad\hence\qquad
(x+1)^3(1-s)=2
$$
(these substitutions are verified directly by use of Lemma~\ref{L1}).
Computing a multiplier of linear transformation between
$\psi$-functions, we obtain finally that functions
$$
\psi_{1,2}^{}=\frac{\sqrt{(x+1)^3-1}}{x+1}\left\{
P_{\!\!\sm\frac13}\!\!\!\!\left(1-\Mfrac{2}{(x+1)^3} \right)\!,
\quad
Q_{\!\sm\frac13}\!\!\!\!\left(1-\Mfrac{2}{(x+1)^3} \right)
\right\}
$$
provide a basis of solutions to equation \eqref{2'}.
\end{pf}

We conclude this section with one example which will be used in the
last section (Sect.~\ref{schwarz}) when appearing a hyperelliptic
curve.

\begin{example}\label{C1}
Rational parametrization of the zero genus equivalence \eqref{zero}
generates, as before in Examples~\ref{E2} and \ref{E3}, Fuchsian
equations for uniformizer $\T$. We obtain here the nice equation
\begin{equation}\label{TT}
\psi''=-\frac14\frac{(\T^6-20\,\T^3-8)^2}
{\T^2(\T^3+8)^2(\T^3-1)^2}\,\psi
\end{equation}
defining  monodromy of an 8-punctured sphere; the equation comes from
the change $\T^3=x$ in Chudnovsky's equation $\eqref{3'}$. That this
$x$ is a perfect cube means that Hauptmodul $\T=\T(\tau)$ is certain to
have an explicit representation in terms of classical $\vartheta$- or
Dedekind's eta-functions. This is so indeed and using some results of
work \cite{maier}, one can derive that
\begin{equation}\label{T}
\T(\tau)=-2\,\frac{\eta^3(2\tau)}{\eta^3(\tau)}
\frac{\eta(3\tau)}{\eta(6\tau)}\,,
\end{equation}
where $\eta(\tau)\DEF\prod_k(1-\re^{\ri k\tau})$,
$\tau\in\Hp$\label{page}. It follows that this $\T(\tau)$ satisfies the
equation
$$
\frac{\{\T,\tau\}}{\dot\T{}^2}=-\frac12\frac{(\T^6-20\,\T^3-8)^2}
{\T^2(\T^3+8)^2(\T^3-1)^2}\,.
$$
\end{example}

\section{Chudnovsky's equations and punctured tori\label{s3}}
\subsection{Equations on tori}
Recall that differential equation  on torus is,
by definition, an ODE of the (normal) form
$$
\psi_{\u\u}^{}=\Xi(\u)\,\psi
$$
with some function $\Xi(\u)$ being an elliptic (transcendental) one in
variable $\u$. In order this equation be of Fuchsian class, the
$\Xi(\u)$ must have second order poles at most. Hence this equation
should be representable in form of a sum over poles $\u=\alpha$ of
$\Xi(\u)$:
\begin{equation}\label{Xi}
\psi_{\u\u}^{}=\left\{
{\sum}_\alpha\!\big( A_\alpha\,\wp(\u-\alpha)+
C_\alpha\,\zeta(\u-\alpha)\big)+A_{{\s0}}\right\}\psi\,,\qquad
{\sum}_\alpha C_\alpha=0\,,
\end{equation}
where $\wp$ and $\zeta$ constitute, together with the
$\sigma$-function, the standard Weierstrassian basis of the elliptic
theory \cite{WW,halphen,weber} for  equation
\begin{equation}\label{elliptic}
\begin{aligned}
\wpp^2&=4\,\wp^3-a\,\wp-b\\
&=4\,(\wp-e)(\wp-\ep)(\wp-\epp)\,.
\end{aligned}
\end{equation}

Along with the preceding sections, we are interested in equations
\eqref{Xi} having a Fuchsian monodromy of finite genus. The case of
non-punctured tori is the `well-trod' domain (the theory of elliptic
functions \cite{WW}) so the torus will be considered to have at least
one puncture; one of the coefficients $A_\alpha$ is equal to
$-\frac14$. The simplest such model is the  singly punctured torus
considered for the first time in the classical work \cite{keen}:
\begin{equation}\label{torus1}
\psi''=-\frac14\big\{\wp(\u;a,b)+A\big\}\,\psi\,.
\end{equation}
On the other hand, finiteness of genus tells us that Riemann surface,
whose fundamental group representation is the monodromy $\G_\u$ to
equation \eqref{Xi}, is always related to a certain \textit{finitely
sheeted} cover $\u\mapsto s$ over this torus (or cover $s\mapsto\u$ by
several copies of the torus). This means that there always exists an
equation
\begin{equation}\label{u}
\Phi(s,\u)\DEF G\big(s,\wp(\u),\wpp(\u)\big)=0
\end{equation}
being polynomial in its $s$, $\wp$, $\wpp$-arguments and realizing this
cover. It is a polynomial in $s$-argument and transcendental function
in $\u$-variable. This is a generic form of covers $\u\rightleftarrows
s$ being an analog of the standard models of algebraic curves given by
polynomial dependencies $P(x,s)=0$.

The minimal possible number of $\u$-sheets branching over $(s)$-plane
is equal to 2 since two is the minimal order of elliptic function.
Therefore simplest covers by tori are the 2-sheeted ones and, hence,
the simplest reduction of \eqref{u} is $R(s)=\wp(\u)$, where $R(s)$ is
any rational function. By the implicit function theorem, branch points
$(\u_k,s_k)$ of the map $s\mapsto \u$ are solutions of equations
$\{\Phi=0,\;\Phi_\u=0\}$ plus separate analysis of the point
$\wp(\u)=\infty$. Hence, the high order rational functions $R(s)$ lead
to a large number of such points and the simplest of the cases is thus
\begin{equation}\label{zu}
s=\wp(\u)\,.
\end{equation}
We may consider this equality as a change $\u\mapsto s$ in
Eq.~\eqref{torus1}. Then it becomes a particular case of the well-known
algebraic form to the famous Lam\'e equation
\cite[Sect.~23$\boldsymbol{\cdot}$4]{WW}
\begin{equation}\label{L2}
\psi''=-\frac{3}{16}\!\left\{\,\frac{1}{(s-e)^2}
+\frac{1}{(s-\ep)^2}+\frac{1}{(s-\epp)^2}-\frac13\,
\frac{5\,s-A}{(s-e)(s-\ep)(s-\epp)}
\right\}\psi
\end{equation}
having the signature $(2,2,2,\infty)$ and its single puncture is
located at point $s=\infty$.

\subsection{On Halphen's transformation}
Halphen \cite{halphen} used further the original trick
$$
s\mapsto x:\qquad\{s=\wp(\u)\}\dashrightarrow\{\u=2\,u\}
\dashrightarrow \{\wp(u)=x\}
$$
to convert equation \eqref{L2} into the form
\begin{equation}\label{P1}
\psi''=-\frac14\!\left\{\,\frac{1}{(x-e)^2}
+\frac{1}{(x-\ep)^2}+\frac{1}{(x-\epp)^2}-
\frac{2\,x-A}{(x-e)(x-\ep)(x-\epp)}
\right\}\psi
\end{equation}
which is our case because \eqref{P1} has the signature
$(\infty,\infty,\infty,\infty)$. It is known that inverse Halphen's
transformation, once applied to algebraic form \eqref{P}, turns it into
equation
$$
p\,\Psi''+\frac12\,p'\,\Psi'+\frac{1}{16}\,(s+\tilde A)\,\Psi=0\,,
\qquad
p\DEF(s-e)(s-\ep)(s-\epp)\,,
$$
whose normal form is \eqref{L2} after a simple adjustment of
parameters.

All this material is classical  \cite[p.~471]{halphen},
\cite[\S\,37]{poole} \cite[p.~185]{chud2}, however, exact correlation
between Lame's equations mentioned above and Eq.~\eqref{torus1}
requires more accurate description. It should be noted also some
ambiguity in work \cite{chud2} which mentions an equivalence between
four punctured sphere \eqref{P1} and 1-punctured torus \eqref{torus1},
whereas their monodromies $\G_\u$ and $\G_x$ have even different ranks;
2 and 3 respectively.

The cover \eqref{u} is never single-sheeted one $s\mapsto \u$.
Therefore group $\G_u$ will be either subgroup of $\G_s$ (\eg, the case
\eqref{zu}\footnote{This exhibits, incidentally, an interesting fact: a
non-free rank 3 group $\G_s$ has a free subgroup $\G_\mathfrak{u}$ of a
smaller rank. Genus of $\G_\mathfrak{u}$ is however not zero but
unity.}) or commensurable with it (general case \eqref{u}). This means
in particular that if we have a correct $A$-parameter for punctured
torus \eqref{torus1}, \ie, $\u(\tau)$ is single-valued, then the map
$\u\mapsto s$ of the form \eqref{zu} yields a single-valued function
$s=\wp\big(\u(\tau)\big)$. We thus obtain a `good' $A$-parameter for
$s$-equation \eqref{L2} from that of $\u$-equation \eqref{torus1}; so
$\G_s$ is a correct monodromy for \eqref{L2}, whereas in the opposite
direction $s\mapsto \u$ we have a $(1\!\mapsto\!2)$-map. As for the
general cover \eqref{u}, both of the maps $s\rightleftarrows\u$ are
always non-single-valued and mutual equivalence of the $A$-parameter
problems for \eqref{L2}, \eqref{P1}, and \eqref{torus1} is not obvious
a priori. Below is a complete and precise formulation.

\begin{theorem}\label{T3}
Halphen's transformation is a transcendental version of birational
transformations between representations of elliptic curves
\eqref{elliptic} in form of covers \eqref{u}. This entails an
equivalence of the $A$-parameter problems for equations \eqref{L2},
\eqref{P1}, and \eqref{torus1} and computability of their
$A$-parameters one through another. The quantities $x$, $s$, and $\u$
as functions of the ratio $\tau=\psi_2^{}/\psi_1^{}$ are single-valued
and computable if one of these functions has been known.
\end{theorem}

\begin{pf}
Let us use the duplication formula for Weierstrass $\wp$-function in
order to treat the Halphen formulae above as the bi-single-valued
(transcendental) transformations between  two models $\Phi_1(x,\u)=0$
and $\Phi_2(s,u)=0$ of the one elliptic curve \eqref{elliptic}:
$$
\Phi_1\!:\; x=\wp\mbig[6](\Mfrac12\u\mbig[6])\,,\qquad
\Phi_2\!:\; s=\wp(2\,u)\,.
$$
Indeed, the equality
$$
\wp(2\,u)=-2\,\wp(u)+\frac{1}{16}\,
\frac{\big(12\,\wp^2(u)-a\big)^2}{\wpp(u)^2}
$$
entails the following single-valued transitions
$(x,\u)\rightleftarrows(s,u)$:
\begin{alignat}{4}
x&=\wp(u)\,,\qquad\qquad &s&=\frac{1}{16}\,\frac{(4\,x^2
+a)^2+32\,b\,x}{4\,x^3-a\,x-b}\,,\label{bi1}\\[1ex]
\u&=2\,u\,,\qquad &u&=\frac12\,\u\,.\label{bi2}
\end{alignat}
These, by Lemma~\ref{L1}, realize explicitly transformations between
equations \eqref{L2}, \eqref{P1}, and \eqref{torus1}. Although function
$x$ is an algebraical one of $s$ it is transcendently single-valued of
the pair $(s,u)$. Owing to  isomorphism \eqref{bi1}--\eqref{bi2}, all
the monodromies $\{\G_x$, $\G_s\}$ are the correct Fuchsian ones of
genus 0 and $\{\G_\u$, $\G_u\}$ are of genus 1 as soon as one of them
has been known to be a correct Fuchsian monodromy. From \eqref{bi1} it
also follows that the free group $\G_x$ is an index 4 subgroup of
non-free group $\G_s$. In a more explicit manner, the proof uses
`Puiseux developments' for $\u=\u(x)$ about points $x=\{e,\ep,\epp\}$.
Inverting the standard series for $\wp$-function \cite{WW}, we get a
series of the type
$$
\frac12\,\u_{\s\pm}=\omega\pm\!\!\!\!
\sqrt[\uproot{1}\sm2]{\!12\,e^2-a\,}\cdot
\left\{2-\frac{4\,e}{12\,e^2-a}\,(x-e)+\cdots\right\}\sqrt{x-e\,}
$$
and the similar series for $u=u_{\s\pm}(s)$. In both of these cases the
square root $\sqrt{x-e\,}$ is represented by a single-valued function
of $\tau$ because $x(\tau)-e$ has an exponential behavior in $\tau$
(due to puncture). In turn, $s(\tau)-e$ is a perfect square
$$
s-e=\mbig[8]\{\frac{(x-e)^2-(e-\ep)(e-\epp)}{2\,\wpp(u)}\,\mbig[8]\}^2
$$
and $s(\tau)$ is an exponent again in the vicinity of $s=\infty$. So
$\u(\tau)$ and $u(\tau)$ are additively automorphic single-valued
functions of $\tau$  (Abelian integrals) and $x(\tau)$, $s(\tau)$ are
purely automorphic single-valued ones. All of them are computable
through any other one by means  of Halphen's transformation itself,
that is \eqref{bi1}--\eqref{bi2}.
\end{pf}

\begin{remark}\label{R5}
From uniqueness of Chudnovsky's list it immediately follows the
uniqueness of the four Lam\'e equations \eqref{L2} of signature
$(2,2,2,\infty)$. Correlating substitutions \eqref{JJ} with
\eqref{bi1}, one can show that all the transitions between equations
\eqref{L2} and \eqref{J} are given by the certain \textit{zero genera}
transformations $F(s,J)=0$. We may of course drop these intermediate
Lam\'e equations and then Halphen's transformation becomes just a
single-valued transition from the torus coordinate $\u$ to the
4-punctured one $x$ by the formula $x=\wp\big(\frac12\u\big)$; this is
checked directly by Lemma~\ref{L1}.
\end{remark}

\subsection{Chudnovsky's equations on tori\label{3.3}}
In view of exclusive character of Eqs.~\eqref{1}--\eqref{4}, it is
useful to display the complete list of associated Fuchsian equations on
tori in an explicit form including their equivalences between each
other. The first two cases are simple and related to equations
\eqref{1}, \eqref{2}; they were obtained in work \cite{keen} based on
some symmetry properties. These cases are equations of the form
\eqref{torus1} with a zero value of the parameter $A$. The two other
ones (most nontrivial) do not appear in any modern reference.

Since parameters $(a,b)$ and singular points of equations
\eqref{1}--\eqref{4}, \eqref{P1} (and consequently the $A$-parameter in
\eqref{P1}) are not invariant quantities, we pass from Weierstrass'
$\wp(\u;a,b)$-representation to the invariant object $\wp(\bu|\mu)$
defined by unique modulus $\mu$. The rule reads as follows
$$
\wp(\u;a,b)=\wp(\u|\om,\omp)\FED\frac{1}{\om^2}\,
\wp(\bu|\mu),\qquad
\bu\DEF\frac{\u}{\om}\,,
$$
where $\mu$ and half-periods $\om$, $\omp$ are computed through the
standard elliptic modular inversion problem. In generic case its
solution is defined by the chain of equations \cite{halphen,WW}
\begin{equation}\label{mod}
J(\mu)=\frac{a^3}{a^3-27\,b^2}\,,\qquad
\om=\pm\sqrt{\frac ab\frac{\g3(\mu)}{\g2(\mu)}}\,,
\qquad\omp=\mu\,\om
\end{equation}
and Weierstrass' modular forms $\g2(\mu)$ and $\g3(\mu)$ have numerous
computational formulae. Most convenient of them are representations in
terms of theta-constants. If we introduce the second Jacobi's constant
$$
\vartheta_2(\tau)\DEF\sideset{}{_k}\sum_{-\infty}^{+\infty}
\re_{\mathstrut}^{\left(k+\frac12\right)^2\pi\ri\tau}
$$
then  one can use the following expressions for these forms
\cite{halphen,weber}:
\begin{equation*}\label{g23}
\begin{aligned}
\g2(\mu)&=\frac{\pi^4}{12}\, \big\{\vartheta_2^8(\mu)+
\vartheta_3^8(\mu)-
\vartheta_2^4(\mu)\,\vartheta_3^4(\mu)\big\}\,,\\[1ex]
\g3(\mu)&=\frac{\pi^6}{432}\,
\big\{\vartheta_2^4(\mu)+\vartheta_3^4(\mu)\big\}
\big\{2\,\vartheta_3^4(\mu)-\vartheta_2^4(\mu)\big\}
\big\{\vartheta_3^4(\mu)-2\,\vartheta_2^4(\mu)\big\}\,.
\end{aligned}
\end{equation*}

In order to derive equations on tori we shift singularities of
Eqs.~$\eqref{1'}$--$\eqref{4'}$ into the Weierstrass form \eqref{P1}
with $e+\ep+\epp=0$ and then compute corresponding Klein's
$J$-invariants. One arrives at four tori with moduli
$\{\ri,\varepsilon,\varrho,\varkappa\}$ \cite{chud2}:
$$
J(\ri)=1\,,\qquad J(\varepsilon)=0\,,\qquad
J(\varrho)=\frac{73^3}{2^4 3^7}\,,\qquad
J(\varkappa)=\frac{2^8 31^3}{3^3 5^3}\,.
$$

\begin{theorem}\label{T4}
Suppose parameters $(a,b,A)$ correspond to equations of the form
\eqref{torus1}. Then the following equations
\begin{alignat*}{4}
J(\ri)\,,\;(a,b,A)&=(4,0,0)\,,&\quad\psi_{\boldsymbol{\u\u}}^{}&=
-\wp(2\bu|\ri)\,\psi\,,
\tag{$1''$}\label{1''}\\[0.5ex]
J(\varepsilon)\,,\;(a,b,A)&=(0,4,0)\,,&\psi_{\boldsymbol{\u\u}}^{}&=
-\wp(2\bu|\varepsilon)\,\psi\,,
\tag{$2''$}\label{2''}\\[0.5ex]
J(\varrho)\,,\;(a,b,A)&=\!\left(\frac{292}{3},-
\frac{4760}{27},\frac23\right)\!,
&\psi_{\boldsymbol{\u\u}}^{}&= -\!\bigg\{\wp(2\bu|\varrho)-
\frac16\,\pi^2\,\vartheta_2^4(\varrho)\bigg\}\,\psi\,, \tag{$3''$}
\label{3''}\\[0.5ex]
J(\varkappa)\,,\;(a,b,A)&=\!\left(\frac{496}{3},-\frac{11044}{27},
\frac43\right)\!,\quad &\psi_{\boldsymbol{\u\u}}^{}&=
-\!\bigg\{\wp(2\bu|\varkappa)- \frac{\sqrt{5\,}}{75}\,\pi^2
\,\vartheta_3^4(\varkappa) \bigg\}\,\psi \tag{$4''$}\label{4''}
\end{alignat*}
are the complete set of Fuchsian equations on tori being pullback of a
$_2F_1$-equation by rational functions of $x$; the intermediate
Halphen's transformation $x=\omega^{\sm2}\wp(\bu|\omp\!/\omega)$ is
assumed to be applied.
\end{theorem}

\begin{pf}
Clearly, only two last equations need to be proved. Performing in
\eqref{P1} Halphen's transformation $\om^2x=\wp(\bu|\mu)$, we impart to
Eq.~ \eqref{P1} the form
\begin{equation}\label{Am}
\psi''=-\big\{\wp(2\,\bu|\mu)+\om^2 A(\mu)\big\}\,\psi
\end{equation}
because $\{\wp(z),z\}=-6\,\wp(2\,z)$. If Weierstrass' roots
$(e,\,\ep,\,\epp)$ and their ordering are known, which is our case,
then standard formulae of the elliptic theory \cite{weber,halphen}
$$
\vartheta_2^4(\mu)=\frac{4}{\pi^2}\,(\epp-\ep)\,\om^2\,,\qquad
\vartheta_3^4(\mu)=\frac{4}{\pi^2}\,(e-\ep)\,\om^2\,,\qquad
\vartheta_4^4(\mu)=\frac{4}{\pi^2}\,(e-\epp)\,\om^2
$$
give linear relations between any pair of $\vartheta$-constants and
values of the $\om$-constant for each  case without resorting to
rooting of a generic $g_{2,3}^{}$-ratio in \eqref{mod}. We find that
\begin{alignat*}{4}
\om=\frac12\,\pi\,\ri\,\vartheta_2^2(\varrho)\quad \mbox{for
$J(\varrho)$}\,, \qquad\qquad
\om=\frac{\sqrt[\leftroot{1}\uproot{1}4]{5\,}}{10}\,\pi\,\ri\,
\vartheta_3^2(\varkappa) \quad \mbox{for $J(\varkappa)$} \,.
\end{alignat*}
Substituting this into \eqref{Am}, we get
Eqs.~\eqref{3''}--\eqref{4''}. Completeness of the list follows from a
completeness of the Beukers--Zagier list \cite[p.~427--428]{beukers}.
\end{pf}

It is interesting to notice that Table~\ref{capt} contains an elliptic
curve that does not appear in this theorem; this is the  curve
\eqref{J123}. What is its relation to these tori? To answer this
question let us consider equation \eqref{J123} and derive Fuchsian
equation on torus defined by this curve. It will suffice to use  any of
$x,\,z$-parametrizations of \eqref{J123}:
$$
z=\frac13\,\frac{(3\,\wp(\u)-5)^2}{3\,\wp(\u)+1}\,,\qquad
x=\frac{8}{\wpp(\u)}\,\frac{3\,\wp(\u)-2}{3\,\wp(\u)-5}\,,
$$
where $\wp(\u)\DEF\wp\big(\u;\frac{52}{3},-\frac{280}{27}\big)$, that
is $\wpp^2=\frac{4}{27}(3\,\wp+7)(3\,\wp-5)(3\,\wp-2)$. Applying
Lemma~\ref{L1} with this $z$-change to equation $\eqref{3'}$ (or this
\mbox{$x$-change} to $\eqref{1'}$), we obtain the Fuchsian equation
(changing $\u\mapsto \u-\omega\hbox{\footnotesize$''$}$)
\begin{equation*}\label{p2}
\psi_{\u\u}^{}=-\frac14\!\mbig[9]\{\wp(\u)+
\frac{4}{\wp(\u)-\frac53}+\frac{4}{3} \mbig[9]\}\psi
\end{equation*}
belonging to the general class \eqref{Xi}. This equation has two
punctures at points $\{0,\,\omp\}$ since $\frac53=\wp(\omp)$:
$$
\psi_{\u\u}^{}=-\frac14\!\left\{\wp(\u)+\wp(\u-\omp)-
\frac13\right\}\psi\,.
$$
It therefore reduces to an equation with one puncture  if we make use
of formula for division of the half-period $\omp$ by 2:
$$
\wp(\bu|\mu)+\wp(\bu-\mu|\mu)
=
\wp\big(\bu\big|\tfrac12\mu\big)+\wp(\mu|\mu)\,.
$$
A simple calculation shows that  modulus $\mu$ of this torus is found
to be $\mu=2\,\varrho$; thus, the curve \eqref{J123} does not produce
new $\u$-equation.

\section{Transcendental equivalence\label{trans}}

\subsection{Mutual covers of tori. Examples\label{51}} Just as
equations \eqref{1}--\eqref{4} are equivalent by algebraic
transformations, so are equivalent equations \eqref{1''}--\eqref{4''}.
Their equivalence will be realized by  transcendental changes
$\Xi(\bu,\bs)=0$ coming from Theorem~\ref{T2} and Halphen's
transformations. These changes constitute  mutual covers of tori
$(\bu)$ and $(\bs)$ by each other and are very rich in consequences.
Because of this, we shall not build the `transcendental' analog of
Table~\ref{capt} but restrict ourselves to the most interesting
branches of the previous machinery. In order to exhibit the way of
getting formulae we consider only two exceptional cases of
Theorem~\ref{T2} and, since examples that follow are the first ones
along these lines, expound one of them at greater length.

\begin{example}\label{EE}
As a first instance we derive the transcendental equivalence
$\eqref{1''}\rightleftarrows\eqref{3''}$. Let us start from the
rational (zero genus) counterparts  to Eqs.~$\eqref{1''}_x$ and
$\eqref{3''}_z$. We may perform Halphen's transformations $x\mapsto \u$
in \eqref{1'} and $z\mapsto \bs$ in \eqref{3'} and arrive at a couple
of Fuchsian equations on tori $(\bu)$ and $(\bs)$ whose monodromies, by
virtue of Theorem~\ref{T3}, are known to be Fuchsian. Thus, we put
\begin{equation}\label{xz}
x=\wp(\u|4,0)\FED\frac{1}{\omega^2}\wp(\bu|\ri),\qquad
z+\frac73=\wp\mbig[7](\mathfrak{s}\,\Big
|\Mfrac{292}{3},-\Mfrac{4760}{27}\mbig[7])
\FED\frac{1}{\tilde{\omega}^2}
\wp(\bs|\varrho)\,,
\end{equation}
where constants $\om$ and $\tilde\om$ are the $\om$-constants for
invariants $J(\ri)$ and $J(\varrho)$ respectively. The second of these
tori is, perhaps, not among the exact solvable modular inversion
problems\footnote{We were unable to find out the value $\varrho$ in
tables on imaginary quadratic fields; see, \eg, \cite{weber}.}: we
compute $-\ri\,\varrho\approx 1.563\,\,401\,\,922\!\ldots$ and
$\ri\,\tilde{\omega}\approx 0.539\,\,128\,\,911\!\ldots$. First torus
$\wpp^2=4\,\wp^3-4\,\wp$ is isomorphic to the classical lemniscate
$y^2=x^4-1$ and its $\omega$-constant (the lemniscatic constant) was
obtained by Gauss. In a $\vartheta$-notation, under normalization
$(a,b)=(4,0)$, the constant has the form
$\omega=\frac12\,\pi\,\vartheta_2^2(\ri)\approx
1.311\,\,028\,\,777\!\ldots$.

Now, we consider an algebraic equivalence
$\eqref{1'}\rightleftarrows\eqref{3'}$ determined, say, by formula
\eqref{J123}. Substituting there
$$
x=\om^{\sm2}\wp(\bu|\ri)\,,\qquad
z=\tilde\om^{\sm2}\wp(\bs|\varrho)-\frac{7}{3}\,,
$$
we get
$$
\om^{\sm4}\wp^2(\bu|\ri)\big(1-\om^{\sm4}\wp^2(\bu|\ri)\big)=
\frac{432\,\big(3\,\tilde\om^{\sm2}\wp(\bs|\varrho)-10 \big)}
{\big(3\,\tilde\om^{\sm2}\wp(\bs|\varrho)-7\big)^3
\big(3\,\tilde\om^{\sm2}\wp(\bs|\varrho)+17\big)}
$$
and, since
\begin{align*}
\wpp^2(\bu|\ri)&=4\,
\big(\wp^2(\bu|\ri)-\om^4 \big)\,\wp(\bu|\ri)\,,\\[1ex]
\wpp^2(\bs|\varrho)&=\frac{4}{27}\,
\big(3\,\wp(\bs|\varrho)-7\,\tilde\om^2\big)
\big(3\,\wp(\bs|\varrho)-10\,\tilde\om^2\big)
\big(3\,\wp(\bs|\varrho)+17\,\tilde\om^2\big)\,,
\end{align*}
one derives that
\begin{equation*}\label{13}
-\wp(\bu|\ri)\wpp(\bu|\ri)^2=
\pi^8\vartheta_2^{16}(\ri)\,\tilde{\om}^6\!\left\{
\frac{3\,\wp(\bs|\varrho)-10\,\tilde{\om}^2}
{(3\,\wp(\bs|\varrho)-7\,\tilde{\om}^2)
\wpp(\bs|\varrho)}\right\}^{\!\!2}.
\end{equation*}
We know that Weierstrass' $\wp$-function is a quadratic ratio of
Jacobi's theta-functions plus a branch point $e$ \cite{WW}. Since $e=0$
is one of the branch points for lemniscate, the $\wp(\bu|\ri)$-function
on the left hand side of last equation is a perfect square and,
therefore, the equation itself is reducible. A simple calculation with
theta-functions shows that
$$
\pm\sqrt{\!\wp(\bu|\ri)}=\frac{1}{2}\,\pi\,\vartheta_2^2(\ri)\,
\frac{\theta_3\!\!\!\left(\frac12\bu|\ri\right)}
{\theta_1\!\!\!\left(\frac12\bu|\ri\right)}
$$
under the standard notation \cite{weber,halphen}
\begin{alignat*}{6}
\theta_1(\u|\mu)&\DEF -\ri\,
\sideset{}{_k}\sum\limits_{-\infty}^{+\infty}\! (-1)^k\,
\re^{\left(k+\frac12\right)^2\pi\ri\,\mu}_{\mathstrut}\, \re^{(2k+1)\pi
\ri\, \u}\,,&\qquad& \theta_3(\u|\mu)\DEF
\sideset{}{_k}\sum\limits_{-\infty}^{+\infty}\!\!
\re^{k^2\pi\ri\,\mu}\re^{2k\pi \ri\,\u}\,,\\[1ex]
\theta_2(\u|\mu)&\DEF \phantom{-\ri}\,
\sideset{}{_k}\sum\limits_{-\infty}^{+\infty}\!\!
\re^{\left(k+\frac12\right)^2\pi\ri\,\mu}_{\mathstrut}\re^{(2k+1)\pi
\ri\,\u}\,, && \theta_4(\u|\mu)\DEF
\sideset{}{_k}\sum\limits_{-\infty}^{+\infty}\! (-1)^k\,
\re^{k^2\pi\ri\,\mu}\, \re^{2k\pi \ri\, \u}\,.
\end{alignat*}
Finally, we obtain the sought-for transcendental equivalence of two
(`very simple') equations $\eqref{1''}\rightleftarrows \eqref{3''}$.

\begin{proposition}
The linear transformation $\Psi=\sqrt{\mfrac{d\bu}{d\bs}}\,\psi$ and
finitely-sheeted mutual cover of  tori $(\bu)$ and $(\bs)$$:$
\begin{equation}\label{final}
\Xi(\bu,\bs):\quad
\pm\,2\,\frac{\theta_3\!\!\!\left(\frac12\bu|\ri\right)}
{\theta_1\!\!\!\left(\frac12\bu|\ri\right)}\,\wpp(\bu|\ri)\,
\wpp(\bs|\varrho)=
\pi^6\vartheta_2^6(\ri)\,\vartheta_2^6(\varrho)
\frac{6\,\wp(\bs|\varrho)+5\,\pi^2\vartheta_2^4(\varrho)}
{12\,\wp(\bs|\varrho)+7\,\pi^2\vartheta_2^4(\varrho)}
\end{equation}
$($transcendental change\/$)$ transform Fuchsian equations
\begin{equation*}\label{two}
\Psi_{\!\bu\bu}=
-\wp(2\bu|\ri)\,\Psi\,,\qquad
\psi_{\!\bs\bs}^{}= -\!\left\{\wp(2\bs|\varrho)-
\frac16\,\pi^2
\vartheta_2^4(\varrho)\right\}\psi
\end{equation*}
into each other.
\end{proposition}

Direct check of this statement is a highly nontrivial exercise even
with use of Lemma~\ref{L1}.
\end{example}

\begin{remark}
Transcendental equivalence \eqref{final} differs  from the general one
given by formula \eqref{abcd} because it does not involve
${}_2F_1$-series appearing in $\Psi, \psi$-solutions.
\end{remark}

\begin{example}\label{e4}
We choose an equivalence of Eqs.~$\eqref{2''}_{\bu}$ and
$\eqref{3''}_{\bs}$ defined by the simplest relation of genus zero,
that is \eqref{zero}. In this case we have the equi-anharmonic torus
$\wpp^2=4\,\wp^3-4$. Its $\om$-constant and relation between
$\vartheta$-constants read as follows
$$
\om=\frac16\,
\sqrt[\uproot{1}\leftroot{1}4]{\!-27\,}\,\pi\,
\vartheta_2^2(\varepsilon)\,,\qquad
\vartheta_3(\varepsilon)=
\sqrt[\uproot{1}\leftroot{1}6]{\ri\,}\,\vartheta_2(\varepsilon)
$$
($\om\approx1.214\,\,325\,\,323\!\ldots$). Applying the same technique
as in the previous example, we derive, after a little algebra, one of
the equivalences $\eqref{2''}_\bu\rightleftarrows \eqref{3''}_\bs$:
\begin{equation}\label{99}
\pm\!3\sqrt[\leftroot{1}\uproot{1}4]{\!-3\,}\,
\frac{\theta_2\theta_3\theta_4}{\theta_1\theta_1\theta_1}\!\mbig[5](
\mfrac12\bu\big|\varepsilon\mbig[5])=
\frac{\theta_1^2\,\theta_3^{}}{\theta_4^2\,\theta_2^{}}\!\mbig[5](
\mfrac12\bs\big|\varrho\mbig[5])\,,
\end{equation}
(no $\vartheta$-constants here at all) where we performed an additional
simplification by converting all the Weierstrassian functions into
Jacobi's theta's. This result can also be treated as the fact that
relation \eqref{99} represents a finitely-sheeted mutual cover of two
punctured tori whose global coordinate functions $\bu=\bu(\tau)$ and
$\bs=\bs(\tau)$ satisfy the two autonomic ODEs
\begin{equation*}\label{su}
\{\bu,\tau\}\,\dot{\bu}^{\sm2}=-2\,\wp(2\,\bu|\varepsilon)\,,\qquad
\{\bs,\tau\}\,\dot{\bs}^{\sm2}= -2\,\wp(2\,\bs|\varrho)+
\frac13\,\pi^2\,\vartheta_2^4(\varrho)\,.
\end{equation*}
\end{example}

\subsection{Transcendental automorphism and Abelian integral\label{52}}
We have shown above that there are nontrivial \textit{algebraic}
automorphisms of Eqs.~\eqref{1'}--\eqref{4'}. Using Halphen's
transformation and Theorem~\ref{T2}, we deduce that there are
\textit{transcendental}  automorphisms between equations
\eqref{1''}--\eqref{4''}.  Here is one nice example based on the
elliptic curve \eqref{11}. As before, we obtain
\begin{equation}\label{rr}
\wp(\u;4,0)\,\wpp(\u;4,0)^2\cdot
\wp(\mathfrak{s};4,0)\,\wpp(\mathfrak{s};4,0)^2=1\,.
\end{equation}
Further analysis of this example leads a remarkable consequence which
we are about to exhibit below.

By Theorem~\ref{T3} functions $\u=\u(\tau)$ and
$\mathfrak{s}=\mathfrak{s}(\tau)$ satisfy a common nonlinear 3rd order
\ode. Is it possible to get  analytic formulae to its solutions in
terms of known functions?

Inversion $x=\chi(\tau)$ of the ratio \eqref{ratio} for equation
\eqref{1'} is known. This is a square of Legendre's modulus
$k^2(\tau)=\vartheta_2^4(\tau)\big/\!\vartheta_3^4(\tau)$
\cite{WW,weber,halphen}. Insomuch as we deal with  automorphism,  the
second function $z(\tau)$ should be the same as $\chi(\tau)$ with the
difference that its argument is merely a linear fractional function of
the $\tau$-argument for $\chi(\tau)$. Some routine computations with
$\vartheta$-constants show that
\begin{equation}\label{par}
x=\frac{\vartheta_2^2(\tau)}{\vartheta_3^2(\tau)}\,,\qquad
z=\frac{\vartheta_2^2\big(\frac{\tau-1}{\tau+1}\big)}{
\vartheta_3^2\big(\frac{\tau-1}{\tau+1}\big)}\,,
\end{equation}
that is a parametrization of the lemniscate \eqref{11}. The Halphen
transformation $x=\wp(\u;4,0)$ tells us that $\u$ is an everywhere
finite quantity for all $x$:
\begin{equation}\label{ppp}
\pm\u=\int\limits_\infty^{\;x}\!\!\!
\frac{ds}{\sqrt{4\,s^3-4\,s\,}}=\cdots
\end{equation}
So we shall find $\u=\u(\tau)$ if we can represent this integral in
terms of known functions. Changing here  the integration variable
$s\mapsto \!\!\sqrt[\uproot{1}\sm2]{\!s\,}$, we get\footnote{It
follows, incidentally, that the nice identity $\wp\mbig({}_2F_1\!
\big(\frac12,\frac14;\frac54\big|z\big);4\,z,0\mbig)=1$ holds for all
$z$.}
\begin{equation}\label{ind}
\cdots =-\frac14\!\!\!\!\!\int\limits_0^{\;\;1\!/\!x^2}\!\!\!\!\!
s^{\sm\frac34}_{\mathstrut}(1-s)^{\sm\frac12}_{\mathstrut}\,ds=
\!\!\!\sqrt[\uproot{1}\sm2]{x\,}
\cdot{}_2F_1\!
\mbig[6](\Mfrac12,\Mfrac14;\Mfrac54\Big|\Mfrac{1}{x^2}\mbig[6]),
\end{equation}
since \cite[Sect.~2.2.2]{bateman2}
\begin{equation}\label{Fab}
\int\limits_0^{\;z}\!\!s^{\alpha\sm1}(1-s)^{\sm\beta}\,ds=
\frac1\alpha\,z^\alpha\cdot{}_2F_1 (\beta,\alpha;\alpha+1|z)\,,\qquad
\mathrm{Re}\,(\alpha)>0\,.
\end{equation}
An important point here is the fact that this ${}_2F_1$-representation
for indefinite integral \eqref{ind} should be understood as a
complex-valued analytic function being an additive one with respect to
periodicity moduli for integral  \eqref{ppp}. Insomuch as \eqref{ppp}
or \eqref{ind} is an elliptic integral, it has only two independent
periods \cite{WW,weber} and we assign them to integration over segments
$s\in[0,1]$ and $s\in(-\infty,1]$:
$$
\frac14\!\!\int\limits_0^{\;1}\!\!\!
s^{\sm\frac34}_{\mathstrut}(1-s)^{\sm\frac12}_{\mathstrut}\,ds\FED
{\scriptstyle \Pi}\,,
\qquad
\frac14\!\!\!\int\limits_{-\infty}^{\;1}\!\!\!\!
s^{\sm\frac34}_{\mathstrut}(1-s)^{\sm\frac12}_{\mathstrut}\,ds=
-\ri\,{\scriptstyle \Pi}\,.
$$
Moreover, the integral is a lemniscatic  one; this being so, its
periods must be combinations of the  $\om$-constant appearing in
Example~\ref{EE}. This is so indeed and we found that
\cite[Sect.~22$\boldsymbol{\cdot}$8]{WW}
$$
{\scriptstyle \Pi}=
\sqrt{\frac{\pi^3}{8}}\,\,
\Gamma\!\Big(\Mfrac34\Big)^{\!\!\!\!\sm2}
\approx1.311\,\,028\,\,777\!\ldots,
$$
\ie, ${\scriptstyle \Pi}$ coincides with the lemniscatic $\om$, as it
should. Correlating now \eqref{par} and \eqref{ppp}--\eqref{ind} and
gathering all the remaining constants, we obtain, upon simplification,
the following result.

\begin{proposition}\label{pp}
The two additively automorphic functions
\begin{equation}\label{param}
\bu(\tau)=\frac{2}{\pi\,\vartheta_2^2(\ri)}
\frac{\vartheta_3(\tau)}{\vartheta_2(\tau)}\cdot
{}_2F_1\!\!\mbig[7](\Mfrac12,\Mfrac14;\Mfrac54\Big|
\Mfrac{\vartheta_3^4(\tau)}{\vartheta_2^4(\tau)}\mbig[7]),\qquad
\bs(\tau)=\bu
\Big(\Mfrac{\tau-1}{\tau+1}\Big)
\end{equation}
satisfy the common ODE
\begin{equation*}\label{ode}
\{\bu,\tau\}=
-2\,\wp(2\,\bu|\ri)\,\dot{\bu}^2
\end{equation*}
and turn the $\eqref{1''}_{\bu}\rightleftarrows
\eqref{1''}_{\bs}$-automorphism
\begin{equation}\label{us}
\pm\!8\,\frac{\theta_2^{}\theta_3^2\theta_4^{}}{\theta_1^4}
\!\mbig[5](\mfrac12\bu\big|\ri\mbig[5])
=
\frac{\theta_1^4}{\theta_2^{}\theta_3^2\theta_4^{}}
\!\mbig[5](\mfrac12\bs\big|\ri\mbig[5])
\end{equation}
into identity in variable $\tau$.
\end{proposition}

Expression \eqref{us} was obtained, as in Example~\ref{e4}, by a
$\theta$-simplification of reducible equality \eqref{rr}. Complete
verification of this statement   is also a good exercise\footnote{To
ensure against typos  we had tested this proposition numerically.}. It
is worth to be noticed that  expression \eqref{param} is the first
instance of  analytic formula for the additively automorphic object
(Abelian integral) on a Riemann surface---more precisely, on
orbifold---of a negative (non-zero) curvature. To the best of our
knowledge no one explicit formula of such a kind was known hitherto.

Let us say a few words concerning computing genera of  covers
$\Xi(\bu,\bs)=0$ in example of Eq.~\eqref{us}. First of all we
establish the common base periods for a $\theta$-ratio on left hand
side of this equation and derive that they are the same  as for
function $\wp(\bu|\ri)$. Hence variables $\bu$, $\bs$ are assumed to
belong to the square formed by vertices $(0,2,2\,\ri,2+2\,\ri)$ and
\eqref{us} defines a transcendental (4\,:\,4)-cover. We need to
determine its branch points $(\bu_k,\bs_k)$ and their ramification
indices $q_k^{}$. Based on the implicit function theorem, form, as
usual, equations $\Xi(\bu,\bs)=0$ and $\Xi_\bu(\bu,\bs)=0$; their
compatibility condition then gives equations defining these points. We
have a separation of variables $\mathfrak{U}(\bu)=\mathfrak{S}(\bs)$
and this simplifies computation of genus as before in case of pure
algebraic equations. An easy calculation yields
\begin{alignat*}{3}
\theta_1^{}\! \mbig[5](\mfrac12\,\bu_1\big|\ri\mbig[5])&=0\,,&
\theta_3^{}\! \mbig[5](\mfrac12\,\bu_2\big|\ri\mbig[5])&=0\,,\qquad&
\pi^2\,\vartheta_2^4(\ri)\,\theta_3^4\!
\mbig[5](\mfrac12\,\bu_3\big|\ri\mbig[5])=
2\,\theta_1^4\!\mbig[5](\mfrac12\,\bu_3\big|\ri\mbig[5])\,.\\%[1ex]
\bu_1&=0\,,\qquad&\bu_2&=\tau+1\,,\qquad
\end{alignat*}
One has three points over $\bu_1$: $\bs=\{1,\tau,\tau+1 \}$ and the
respective indices $q_k^{}$ are $\{4\}$, $\{2,2\}$, and $\{4\}$. There
are no ramifications over point $\bu_2$ since Eq.~\eqref{us} has a
structure $\theta_3^2\sim\theta_1^4$ in the vicinity of this place. The
four remaining points $\bu_3$, as the local analysis shows, turn out to
be just points of regularity: $q_k^{}=\{1,1,1,1\}$ there. Now, using
the Riemann--Hurwitz formula $g=\frac12\sum(q_k^{}-1)+N(g'-1)+1$ with
$N=4$ and $g'=1$ (torus being covered), we obtain
$$
g=\frac{1}{2}\,\big\{ (2-1)\,2+(4-1)\,2 \big\}+4\,(1-1)+1=5
$$
and arrive again at a Riemann surface of genus five.

Analogs of formulae \eqref{ppp}--\eqref{us} for the equi-anharmonic
case are derived in a similar manner. Lemniscatic  and equi-anharmonic
cases are the only ones  we were able to obtain explicit analytic
formulae.

\subsection{A hyperelliptic curve\label{schwarz}}

The remarkable fact is that  covers and Halphen's transformations
provide  the independent ways of generation of algebraic curves and
integrable Fuchsian equations with finite genus monodromies. All this
is obtained by correlating Chudnovsky's curves $F(x,z)=0$, base covers
of tori $\{x=\wp(\bu), z=\tilde\wp(\bs)\}$, and transcendental covers
$\Xi(\bu,\bs)=0$ between themselves. Lack of space prevents us giving
an exhaustive analysis and we restrict ourselves to considering a
distinguishing example.

\begin{example}
Let us consider transcendental counterpart of curve \eqref{zero}, that
is Eq.~\eqref{99}. Its genus is easily counted because it is seen at
once that ramifications $\bu=\bu(\bs)$ are possible only at place
$\theta_1\big(\frac12\,\bu|\varepsilon\big)=0$. Right hand side of
\eqref{99} tells us that there are only two points over this $\bu$,
namely, points determined by equations
$\theta_4\big(\frac12\,\bs|\varrho\big)=0$ and
$\theta_2\big(\frac12\,\bs|\varrho\big)=0$. Both of their indices are,
obviously, $q=\{3\}$. Riemann--Hurwitz formula above shows, thus, that
genus of \eqref{99}, as a (3\,:\,3)-cover, is $g=3$. What can we say
about algebraic models to this cover?

Replacing variables $x\mapsto x-1$, $z\mapsto z-\frac73$ and
introducing the second coordinates of tori  as
$\wpp(\bu|\varepsilon)\FED y$ and $\wpp(\bs|\varrho)\FED
4\,\sqrt{3\,}\,w$, we may rewrite \eqref{99} as follows:
$$
x^3=\frac{(3\,z-1)^3}{(3\,z+17)(3\,z-10)^2}\,,\qquad
y^2=4\,x^3-4\,,\qquad 324\,w^2=(3\,z-7)(3\,z-10)(3\,z+17)\,.
$$
This 1-dimensional surface in a 4-space $(x,y,z,w)$ contains the plane
$\{ (y,z)$, $(y,w)$, $(x,w)\}$-curves of respective genera
$g=\{0,1,3\}$. Of course, they can not appear in the previous analysis
but we can do birational transformations. Doing that, we observe that
the genus $g=1$ curve $F(y,w)=0$ is isomorphic to the curve
\eqref{J123} with a duplicate modulus $\mu=2\,\varrho$ and the genus
$g=3$ curve
$$
3^{\sm3}\,(x^3-1)^2\big((x^3-1)\,w^2+9\big)w^4=
2^{\sm4}(x^6+64\,x^3+16)\,w^2+x^3-1
$$
is isomorphic to the hyperelliptic form $v^2=(u^3+8)(u^3-1)\,u$. We do
not display here all these birationalities. One immediately recognizes
that the branch $u$-points are  singularities of Eq.~\eqref{TT}.  On
the other hand, Table~\ref{capt} contains two genus 3 curves and,
curiously, we found that one of them is also hyperelliptic; this is the
$\eqref{3'}_z \rightleftarrows\eqref{1'}_x$ curve. Moreover, both of
these hyperelliptic surfaces are isomorphic to one another and can be
turned into the famous classical Schwarz form
$$
\boldsymbol{z}^2=\bx^8+14\,\bx^4+1\,.
$$

The last step is Fuchsian equations. Lemma~\ref{L1}, upon application
of the chain of transitions $x\mapsto w\mapsto u\mapsto \bx$, gives,
however, a Fuchsian $\bx$-equation with singularities located at roots
$(\alpha^8+14\,\alpha^4+1)(\alpha^5-\alpha)=0$ with `excessive' roots
$\alpha^5-\alpha=0$. Nevertheless, we can use of \eqref{TT} because it
has also been arisen from the $\eqref{3'}_z
\rightleftarrows\eqref{1'}_x$-curve \eqref{zero}. Different ways
produce different Fuchsian equations but one curve.

Finally, insomuch as the explicit form for Schwarz's Hauptmodul
$\bx(\tau)$ does not appear in the literature, it is pertinent to
present here its `non-excessive' form:
$$
\bx(\tau)=\frac{1+\ri}{2}\,\frac{(1+\sqrt{3})\,\T(\tau)+2}
{\T(\tau)-1-\sqrt{3}}\,,
$$
where $\T(\tau)$ is defined by formula \eqref{T}. By this means the
transformation $\T\mapsto \bx$ leads to a Fuchsian equation with eight
parabolic singularities $\alpha^8+14\,\alpha^4+1=0$. The equation is
obtained by use of Lemma~\ref{L1} applied to \eqref{TT}:
$$
\psi''=\frac{48\,(\bx^5-\bx)^2}{(\bx^8+14\,\bx^4+1)^2}\,\psi.
$$
\end{example}

\section{Conclusive remarks}
An abundance of Riemann surfaces/orbifolds coming form Chudnovsky's
equations is a nontrivial result in its own right and all this material
requires  development of an independent theory explaining the genesis
of a huge variety of the surfaces, further classification, and with it
unification of getting the formulae. Proposition~\ref{pp} is not
restricted to the elliptic and holomorphic integrals. As it follows
from formula~\eqref{Fab} any holomorphic or meromorphic Abelian
integral
$$
\mathfrak{A}=\int\limits^{\;\,z}\!\!\! s^{\frac kn}_{\mathstrut}
(s^m-1)^{\frac \ell n}_{\mathstrut}\,ds
$$
belonging to the algebraic irrationality $w^n=z^k(z^m-1)^\ell$ with
three branch points $z=\{0,1,\infty\}$ can be worked out in a similar
manner
$$
\mathfrak{A}\dashrightarrow\mathfrak{A}(\tau)\sim {}_2F_1(\chi(\tau))\,,
\qquad\{\mathfrak{A},\tau\}=\Xi(\mathfrak{A})\,\dot{\mathfrak{A}}^2
$$
if Hauptmodul $z=\chi(\tau)$ is known. This is frequently our case;
\eg,  automorphisms considered in Sect.~\ref{aut} lead to curves
$w^n=(z^6-1)^\ell$ and $w^n=z^k\,(z^4-1)^\ell$.

\end{document}